\documentclass[12pt]{amsart}
\usepackage{amsmath,latexsym,amsxtra,amssymb}
\usepackage{bbm}
\allowdisplaybreaks

\newtheorem{lemma}{Lemma}
\newtheorem{prop}[lemma]{Proposition}
\newtheorem{theorem}[lemma]{Theorem}

\newtheorem{remark}[lemma]{Remark}

\newtheorem{exam}[lemma]{Example}

\numberwithin{equation}{section}
\numberwithin{lemma}{section}

\newcommand{\vs}{\left ( }
\newcommand{\os}{\right ) }
\newcommand{\vhaka}{\left [ }
\newcommand{\ohaka}{\right ] }
\newcommand{\vits}{\left |}
\newcommand{\oits}{\right |}
\newcommand{\vnorm}{\left | \left |}
\newcommand{\onorm}{\right | \right |}
\newcommand{\vaalto}{\left \{ }
\newcommand{\oaalto}{\right \} }

\newcommand{\vph}{\varphi}
\newcommand{\Tau}{\mathcal{T}}
\newcommand{\eps}{\varepsilon}

\newcommand{\oA}{{\mathcal{A}}}
\newcommand{\real}{{ \mathbb R}}
\newcommand{\natu}{{ \mathbb N}}
\newcommand{\odotu}{\mathbbm{E}}

\newcommand{\indi}[1]{\mathbbm{1}_{#1}}
\newcommand{\borel}[1]{\mathcal{B}(#1)}

\newcommand{\deri}[2][]{\frac{\partial^{#1}}{\partial #2^{#1}}}
\newcommand{\deeri}[2][]{\frac{\partial^{#1}}{\partial #2}}
\newcommand{\tn}{\mathbb{P}}
\newcommand{\probsp}{(\Omega , {\mathcal F},\tn)}

\newcommand{\saF}{{\mathcal F}}

\def\myproof{\par \vskip0.5cm \noindent {\bf Proof.  }}
\newcommand{\stopproof}{\hfill$\Box$\vspace{.5em}}

\begin{document}
\begin{center}
\textbf{On discrete time hedging in $d$-dimensional option pricing models}

\textsc{M. Hujo}\footnote{Address correspondence to the author at Department of Mathematics and Statistics, University of Kuopio, P.O. Box 1627, FIN-70211 Kuopio, Finland; e-mail mika.hujo@gmail.com}

\textit{Department of Mathematics and Statistics, University of Kuopio, Finland}
\end{center}

\begin{center}
\textsc{Abstract}
\end{center}
We study the approximation of certain stochastic integrals with respect to a $d$-dimen\-sio\-nal diffusion by corresponding stochastic integrals with piece-wise constant integrands i.e. an approximation of the form
$
\sum_{k=1}^d \int_0^T N_s^k dX_s^k \approx  \sum_{k=1}^d \sum_{i=1}^n N_{t_{i-1}}^k (X_{t_i}^k - X_{t_{i-1}}^k).
$ In finance this corresponds to replacing a continuously adjusted portfolio by discretely adjusted one. 
The approximation error is measured with respect to $L^2$ and it is shown that under certain assumptions the approximation rate is $n^{-1/2}$ when one optimizes over deterministic but not necessarily equidistant time-nets $0=t_0 \leq t_1 \leq \cdots \leq t_n = T$.\\

\begin{center}
\end{center}

%%%%%%%%%%%%%%%%%%%%%%%%%%%%%%%%%%%%%%%%%%%%%%%%%%%%%%%%%%%%%%%%%%%%%%%%%

\setcounter{lemma}{0}

%%%%%%%%%%%%%%%%%%%%%%%%%%%%%%%%%%%%%%%%%%%%%%%%%%%%%%%%%%%%%%%%%%%%%%%%%
\vspace*{2em}
\section{Introduction}

Assume a Borel-function $f:\real^d \to \real$, $T > 0$ and a stochastic process $(X_t)_{t\in[0,T]}$ defined as a solution of
\begin{equation}\label{eq:introX}
X_t^i = x_0^i + \int_0^t b_i(X_u) du + \sum_{j=1}^d\int_0^t \sigma_{ij}(X_u)dW_u^j,\ i\in
\vaalto 1,\ldots, d\oaalto,
\end{equation}
where $(W_t)_{t\in [0,T]}$ is a $d$-dimensional Brownian motion and the functions $b$ and $\sigma$ satisfy certain assumptions (cf. Chapter \ref{cha:oletukset}).

Consider the problem that a trader has to hedge, by a self-financing strategy, a European type option with maturity $T>0$, where the pay-off of the option is described by a random variable $f(X_T)$. The perfect hedging strategy is determined by the process $(N_u)_{u\in[0,T]}$ in a stochastic integral representation of $f(X_T)$,
\begin{equation*}
f(X_T) = V_0 + \sum_{k=1}^d \int_0^T N_u^k dX_u^k,
\end{equation*}
where $V_0$ is the initial capital.
In practice the continuous strategy has to be replaced by a discretely adjusted one. This leads to an approximation
$$	
 \sum_{k=1}^d \int_0^T N_u^k dX_u^k \approx  \sum_{k=1}^d \sum_{i=1}^n N_{t_{i-1}}^k (X_{t_i}^k - X_{t_{i-1}}^k),
$$
where $0=t_0 \leq t_1 \leq t_2 \leq \cdots \leq t_n = T$ is a deterministic but not necessarily equidistant time-net. 

We will measure and (to some extent) optimize the error of this approximation in $L^2$.
Our interest lies in the rate of convergence of the approximation, when the approximation error is minimized over all time-nets with at most $n+1$ time-knots. This means that we are interested in the quantity
\begin{equation}\label{eq:intro_virhe}
\inf_{\tau\in \Tau_n} \vnorm \sum_{k=1}^d \int_0^T N_s^k dX_s^k -  \sum_{k=1}^d \sum_{i=1}^m N_{t_{i-1}}^k (X_{t_i}^k - X_{t_{i-1}}^k)\onorm_{L^2}
\end{equation}
as $n$ tends to infinity, where
$$
 \Tau_n := \{(t_i)_{i=0}^m:0=t_0< t_1< \cdots < t_m = T, m\leq n\}.
$$

Let us recall some results from the literature. Among others, the 
1-dimensional case has been considered by Zhang \cite{ZhangR}, Gobet-Temam \cite{Gobet_Temam} and
Geiss \cite{Geiss15}. Geiss considered the approximation problem
for general deterministic nets, which are not necessarily 
equidistant, and a closed form formula for the $L_2$-error was obtained.
Based on this, in \cite{Geiss_Hujo} several classes of examples were given, 
where the optimal rate of convergence $n^{-1/2}$ is attained by general 
deterministic nets (but, in general, not by equidistant ones).
The result from \cite{Geiss15} and \cite{Geiss_Hujo} cannot be straightforward 
extended to the multi-dimensional case because part of the arguments 
from the 1-dimensional case do not seem to apply in the multi-dimensional
situation.

The multi-dimensional case was, for example, studied by Zhang \cite{ZhangR} 
and Temam \cite{Temam} for equidistant nets. For $C^1$-functions with derivatives of polynomial growth, cf. \cite[Proposition 3.1.6 and Corollary 3.3.3]{ZhangR} , Zhang established the rate $n^{-1/2}$. On the 
other side, Temam \cite{Temam} proved the rate $n^{-1/4}$ for the European digital option.

The aim of this paper is to improve the approximation rate 
of the European digital option in the multi-dimensional case
from $n^{-1/4}$ to $n^{-1/2}$ by replacing the equidistant nets 
by general deterministic nets.

The paper is organized as follows: In section \ref{cha:oletukset} we explain the setting we are working with. Section \ref{cha:results} introduces Theorem \ref{the:paatulos} which is our main result. In Theorem \ref{the:paatulos} we show for a certain class of functions $f$, including European digital option, that one gets the $L^2$-approximation rate of $n^{-1/2}$ by optimizing over all deterministic nets of cardinality $n+1$. Our Theorem also allows a drift term in the underlying diffusion process (which is sometimes remarked, but not carried out, in the literature). Section \ref{cha:esim} gives some
examples illustrating Theorem \ref{the:paatulos}.

\vspace*{2em}
\section{Preliminaries}\label{cha:oletukset}
In this chapter we introduce the setting we are working with and recall some known facts that are needed in order to prove our results.

We shall use the standard assumptions from stochastic calculus, i.e. we assume a complete probability space $\probsp$ and, for $T>0$, a right-continuous filtration $(\saF_t)_{t\in[0,T]}$ generated by a standard $d$-dimensional Brownian motion $W=(W_t)_{t\in[0,T]}$ such that $\saF_T=\saF$ and $\saF_0$ contains all null-sets of $\saF$ (cf. \cite{Kar-Shreve1}). 
By $\vnorm x \onorm$ we denote the Euclidean norm of a vector $x \in \real^d$. 
A Borel-function $\vph:B\to \real$ on some set $B\subset \real$ will be extended to $B^d \subset \real^d$ by the notation 
$$
\vph(x) := (\vph(x_1), \vph(x_2), \ldots, \vph(x_d)),\ x\in B^d.
$$
We consider a diffusion
\begin{equation}\label{eq:SDE_X}
 X_t^i = x_0^i + \int_0^t b_i(X_u)du + \sum_{j=1}^d \int_0^t \sigma_{ij}(X_u) dW_u^j,\ i=1, \ldots, d,\ a.s.
\end{equation}
where $x_0\in \real^d$. The process $X$ is obtained through $Y$ given as the unique path-wise continuous solution of (cf. \cite[Corollary 2.2.1 on p. 101]{Nualart})
\begin{equation}\label{eq:SDE_Y}
 Y_t^i = y_0^i + \int_0^t \hat{b}_i(Y_u)du + \sum_{j=1}^d \int_0^t \hat{\sigma}_{ij}(Y_u) dW_u^j,\ i=1, \ldots, d,\ a.s.
\end{equation}
where
\begin{equation}\label{eq:kerroinfunktio_oletus}
 \hat{b}_i(x), \hat{\sigma}_{ij}(x) \in C_b^\infty(\real^d)
\end{equation}
and $\hat{\sigma}\hat{\sigma}^T$, where $(\hat{\sigma}\hat{\sigma}^T)_{ij}(x) = \sum_{k=1}^d \hat{\sigma}_{ik}(x) \hat{\sigma}_{jk}(x)$, is uniformly elliptic i.e. 
$$
 \sum_{i,j=1}^d (\hat{\sigma}\hat{\sigma}^T)_{ij}(x) \xi_i \xi_j \geq \lambda \vnorm \xi \onorm^2,\text{ for all } x, \xi\in \real^d \text{ and some } \lambda > 0.
$$
Under these assumptions the process $Y$ has a transition density $\Gamma$ with appropriate tail estimates (see Theorem \ref{the:siirtymatod} in Appendix).

We consider two cases. The first case 
\begin{itemize}
\item[$(C1)$] $x_0=y_0\in \real^d$, $\hat{b}_i(x):=b_i(x)$, $\hat{\sigma}_{ij}(x) := \sigma_{ij}(x)$, $X_t = Y_t$,
\end{itemize}
is related to the Brownian motion and the second case
\begin{itemize}
\item[$(C2)$] $x_0=e^{y_0} \in (0,\infty)^d$, $\hat{b}_i(y):= \frac{b_i(e^{y})}{e^{y_i}} - \frac{1}{2} \sum_{j=1}^d \hat{\sigma}_{ij}^2(y)$, $\hat{\sigma}_{ij}(y) := \frac{\sigma_{ij}(e^y)}{e^{y_i}}$ and $X_t = e^{Y_t}$,
\end{itemize}
with the convention $e^y = (e^{y_1},\ldots, e^{y_d})$ for $y \in \real^d$, is close to the geometric Brownian motion. In both cases we have
\begin{equation}\label{eq:X_normi_rajoitettu}
\odotu \sup_{t\in[0,T]} \vnorm X_t\onorm^p < \infty
\end{equation}
for any $p > 0$ (cf. \cite[Corollary 2.2.1 on p. 101]{Nualart}). 

To summarize the above, we start with the process $X$ by choosing the matrix $\sigma$ and the vector $b$ such that the matrix $\hat{\sigma}$ and the vector $\hat{b}$ satisfy the required conditions above. In this way we obtain the process $Y$ and deduce properties of the process $X$ from the properties of $Y$.

To handle both of these cases simultaneously, we define functions $Q_i:\real^d \to \real$ for $i=1,\ldots, d$ by
\begin{equation*}
Q_i(x) := \vaalto
  \begin{array}{ll}
    1, & \text{ in case } (C1) \\
    x_i, & \text{ in case } (C2).
  \end{array}\right. 
\end{equation*}
In what follows we assume, for some $q\in [2,\infty)$ and $C > 0$, that 
\begin{equation}\label{eq:f_oletus}
\vits f(x) \oits \leq C\vs 1 + \vnorm x\onorm^q\os,\ x \in E,
\end{equation}
where the $f:E \to \real$ is a Borel-function and the set $E$ is defined by
$$
E := \vaalto
  \begin{array}{ll}
    \real^d, & \text{ in case } (C1) \\
    (0,\infty)^d, & \text{ in case } (C2).
  \end{array}\right. 
$$
Through the function $f$ we define the function $g$ on $\real^d$ by
$$
g(y):=\vaalto \begin{array}{ll}
    f(y), & \text{ in case } (C1) \\
    f(e^y), & \text{ in case } (C2). 
  \end{array}\right. 
$$

Applying Theorem \ref{the:siirtymatod} to the stochastic differential equation
$$
\vaalto
  \begin{array}{ll}
Z_t^i = Z_0^i + \sum_{j=1}^d \int_0^t \hat{\sigma}_{ij}(Z_u) dW_u^j, & \text{ in case } (C1) \\
Z_t^i = Z_0^i - \int_0^t \vs \frac{1}{2} \sum_{j=1}^d \hat{\sigma}_{ij}^2(Z_u)\os du + \sum_{j=1}^d \int_0^t \hat{\sigma}_{ij}(Z_u) dW_u^j, & \text{ in case } (C2)
 \end{array}\right.  
$$
gives a transition density $\Gamma_0$ such that we can define the function $G \in C^\infty([0,T) \times \real^d)$ by
\begin{equation*}\label{eq:G}
 G(t,y):= \int_{\real^d} \Gamma_0(T-t,y,\xi) g(\xi) d\xi,\ 0\leq t < T
\end{equation*}
so that
\begin{equation}\label{eq:G_yhtalo}
\vaalto
  \begin{array}{ll}
\vs \deri{t} + \frac{1}{2}\sum_{k,l=1}^d  \vs \hat{\sigma} \hat{\sigma}^T (y)\os_{kl} \deeri[2]{y_k y_l}\os G(t,y) = 0 & (C1) \\
\vs \deri{t} - \sum_{i=1}^d \vs\frac{1}{2} \sum_{j=1}^d \hat{\sigma}_{ij}^2(y)\os \deri{y_i} + \frac{1}{2}\sum_{k,l=1}^d \vs \hat{\sigma} \hat{\sigma}^T (y)\os_{kl} \deeri[2]{y_k y_l}\os G(t,y) = 0 & (C2).
 \end{array}\right.  
\end{equation}

Now we can define the function $F$ on $[0,T)\times E$ by
$$
F(t,x) :=\vaalto \begin{array}{ll}
    G(t,x), & \text{ in case } (C1) \\
    G(t,\log(x)), & \text{ in case } (C2).
  \end{array}\right. 
$$
Assumption (\ref{eq:f_oletus}) together with Theorem \ref{the:siirtymatod} implies that for $0\leq t \leq T' < T$
\begin{equation}\label{eq:ey_F_1_deri}
\vits Q_i (x) \oits \vits \deri{x_i} F(t,x)\oits \leq C_{d,T'} (1 + \vnorm x \onorm^q),\ x\in E, \ i=1,\ldots, d
\end{equation}
and
\begin{equation}\label{eq:ey_F_2_deri}
\vits Q_i (x) \oits\vits Q_j (x) \oits \vits \deeri[2]{x_i x_j} F(t,x)\oits \leq C_{d,T'} (1 + \vnorm x \onorm^q),\ \ x\in E,\ i,j=1,\ldots, d.
\end{equation}

Let
\begin{equation}\label{eq:oA}
\oA:=\deri{t} + \frac{1}{2}\sum_{k,l=1}^d A_{kl}(x) \deeri[2]{x_k x_l}
\end{equation}
where
\begin{equation}\label{eq:matriisi_A}
A_{kl}(x):=\sum_{j=1}^d \sigma_{kj}(x)\sigma_{lj}(x).
\end{equation}
From the definition of $F$ and equation (\ref{eq:G_yhtalo}) it follows that
\begin{equation}\label{eq:bpde}
\oA F(t,x)=0 \text{ on } [0,T)\times E.
\end{equation}
Moreover, It\^o's formula gives that
$$
F(t,X_t) = F(0,X_0) + \sum_{k=1}^d \int_0^t \deri{x_k} F(u,X_u) dX_u^k,\ a.s. \ t \in [0,T).
$$
Finally, Theorem \ref{the:siirtymatod} gives that
\begin{equation*}
F(t,X_t) \to f(X_T) \text{ in } L^2 \text{ as } t\nearrow T
\end{equation*}
and 
$$
f(X_T) = F(0,X_0) + \sum_{k=1}^d \int_0^T \deri{x_k} F(u,X_u) dX_u^k\ a.s.
$$

\vspace*{2em}
\section{Results}\label{cha:results}

In the rest of the paper we assume the setting from Chapter \ref{cha:oletukset}. We start this chapter by stating our main result Theorem \ref{the:paatulos}. It implies that under certain conditions the convergence rate for the supremum of the approximation error is bounded by $n^{-1/2}$, when one optimizes over all deterministic time-nets of cardinality $n+1$. Two examples where Theorem \ref{the:paatulos} is applied to are presented in Chapter \ref{cha:esim}.

\begin{theorem}\label{the:paatulos}
Assume that for all $x \in E$
$$
\vits \frac{\partial^s}{\partial_{x_\beta}^q \partial_{x_\alpha}^r} \sigma_{ij}(x)\oits \leq C_1 \frac{Q_i(x)}{Q_\beta^q(x) Q_\alpha^r(x)}, \text{ where } q+r = s,\ q,r,s\in \vaalto 0,1,2\oaalto,
$$
$\vits b_i(x)\oits \leq C_1 Q_i(x)$ and $ A_{ii}(x) \geq \frac{1}{C_1}Q_i^2(x)$
for $i\in \vaalto 1,\ldots, d\oaalto$ and some fixed $C_1 > 0$. Moreover, assume that
\begin{equation}\label{eq:tulos2oletus}
\sup_{\alpha,\beta}
\odotu \vhaka A_{\alpha \alpha} (X_t) A_{\beta \beta}(X_t) \vits \deeri[2]{x_\alpha x_\beta} F(t,X_t)\oits^2 \ohaka \leq \frac{C_2}{(T-t)^{2\theta}},\ \theta \in [0,1), \text{ for some } C_2>0.
\end{equation}
Then
$$
\vs \odotu \sup_{t\in[0,T]} \vits \sum_{i=1}^n \sum_{k=1}^d \int_{t_{i-1}^\eta \wedge t}^{t_i^\eta \wedge t} \vs\deri{x_k} F(u,X_u) - \deri{x_k} F(t_{i-1}^\eta,X_{t_{i-1}^\eta})\os dX_u^k \oits^2\os^\frac{1}{2} \leq \frac{D_1}{\sqrt{n}},
$$
where
$$
\tau_n^\eta = (t_i^\eta)_{i=0}^n := \vs T\vs 1 - \vs 1 - \frac{i}{n}\os^\frac{1}{1-\eta} \os\os_{i=0}^n \text{ and }
\vaalto \begin{array}{ll}
    \eta =0, & \theta \in [0,\frac{1}{2}) \\
    \eta\in (2\theta - 1, 1), & \theta \in [\frac{1}{2},1) \\
  \end{array}\right.
$$
and $D_1>0$ depends at most on $\eta, C_1, C_2, d$ and $T$.

In addition, assume that
\begin{equation}\label{eq:oletus_H}
\inf_{u\in (r,s)} H^2(u)=C_H> 0,
\end{equation}
for some $0 \leq r < s < T$, where $H$ is defined by 
\begin{equation}\label{eq:H}
H^2(u):=\odotu \sum_{\alpha, \beta, i, k =1}^d A_{\alpha \beta}(X_u) A_{i k}(X_u) \deeri[2]{x_\alpha x_i}F(u,X_u) \deeri[2]{x_\beta x_k}F(u,X_u), u \in [0,T).
\end{equation}
Then we have the following two cases:
\begin{enumerate}
\item[($L1$)]
In the case that $\theta \in [0,3/4)$, we have, for any sequence of time-nets $0= t_0^n\leq t_1^n \leq \ldots \leq t_n^n = T$ with $\sup_{i=1,\ldots,n} (t_i^n - t_{i-1}^n) \leq C_\tau/n$, $C_\tau >0$, that
\begin{align}\label{eq:alaraja}
&\liminf_{n\to \infty} \sqrt{n} \vs \odotu \sup_{t\in[0,T]} \vits \sum_{i=1}^n \sum_{k=1}^d \int_{t_{i-1}^n\wedge t}^{t_i^n\wedge t} \vs\deri{x_k} F(u,X_u) - \deri{x_k} F(t_{i-1},X_{t_{i-1}})\os dX_u^k \oits^2\os^\frac{1}{2} \geq \frac{1}{D_2}.
\end{align}
\item[($L2$)] 
If $\theta \in [3/4, 1)$, then we have that
\begin{align}\label{eq:alaraja2}
&\liminf_{n\to \infty} \sqrt{n} \vs \odotu \sup_{t\in[0,T]} \vits \sum_{i=1}^n \sum_{k=1}^d \int_{t_{i-1}^{\eta,n}\wedge t}^{t_i^{\eta,n}\wedge t} \vs\deri{x_k} F(u,X_u) - \deri{x_k} F(t_{i-1}^{\eta,n},X_{t_{i-1}^{\eta,n}})\os dX_u^k \oits^2\os^\frac{1}{2}\\\nonumber
&\geq \frac{1}{D_2}.
\end{align}
The constant $D_2>0$ depends at most on $C_1, C_2, C_H, d$ and $T$.
\end{enumerate}
\end{theorem}

\begin{remark}
\begin{enumerate}
\item In the case that the process $(X_t)_{t\in [0,T]}$ does not have a drift, it follows from Doob's inequality that inequalities (\ref{eq:alaraja}) and (\ref{eq:alaraja2}) can be replaced by 
\begin{equation*}
\liminf_{n\to \infty} \sqrt{n} \vs \odotu \vits \sum_{i=1}^n \sum_{k=1}^d \int_{t_{i-1}^n}^{t_i^n} \vs\deri{x_k} F(u,X_u) - \deri{x_k} F(t_{i-1}^n,X_{t_{i-1}^n})\os dX_u^k \oits^2\os^\frac{1}{2}\geq \frac{1}{4 D_2}.
\end{equation*}
Of course in case of (\ref{eq:alaraja2}) we set $t_i^n = t_i^{\eta,n}$
\item In $(L2)$ we have the lower bound only for time-nets $\tau^\eta$. Compared to $(L1)$ this does not seem natural, since larger $\theta$ should correspond to a worse approximation. We need this restriction for technical reason (but believe that it can be removed).
\item Under the setting of the Chapter \ref{cha:oletukset} the assumptions in Theorem \ref{the:paatulos} which concern the estimates of the matrices $A$ and $\sigma$ and the vector $b$ by the functions $Q_i$ are always satisfied for some $C_1 >0$.
\item It follows by a simple calculation that
\begin{equation*}
H^2(u) = \odotu \sum_{m,n=1}^d\vs \sum_{\alpha, \beta = 1}^d \sigma_{\alpha m}(X_u) \sigma_{\beta n}(X_u) \deeri[2]{x_\alpha x_\beta}F(u,X_u)\os^2.
\end{equation*}
Now because of (\ref{eq:ey_F_2_deri}) we have that $H^2(u) \in [0,\infty)$, for $u\in [0,T)$.
\item If the matrix $A$ defined in (\ref{eq:matriisi_A}) is a diagonal matrix, then 
$$
H^2(u) = \odotu \sum_{\alpha, \beta =1}^d A_{\alpha \alpha}(X_u) A_{\beta \beta}(X_u) \vits \deeri[2]{x_\alpha x_\beta}F(u,X_u)\oits^2
$$
and thus it is equivalent to the function
$$
\sup_{\alpha,\beta}
\odotu \vhaka A_{\alpha \alpha} (X_t) A_{\beta \beta}(X_t) \vits \deeri[2]{x_\alpha x_\beta} F(t,X_t)\oits^2 \ohaka
$$
considered for the upper bound in Theorem \ref{the:paatulos}. 
In the 1-dimensional case our function $H$ is the same as the function $H$ controlling the approximation error in \cite{Geiss15}.
\end{enumerate}
\end{remark}

Now turn to the proof of Theorem \ref{the:paatulos}. We deal with a multi-step approximation error i.e. the stochas\-tic integral $\sum_{k=1}^d \int_0^T \deri{x_k} F(u,X_u) dX_u^k$ is approximated by the stochastic integral $\sum_{k=1}^d \sum_{i=1}^n \int_{t_{i-1}}^{t_i} \deri{x_k} F({t_{i-1}},X_{t_{i-1}}) dX_u^k$. In order to estimate the multi-step error we need to have information about the one-step error occurring in a time interval $[t_{i-1}, t_i]$. Here Proposition \ref{pro:paatulos} and Proposition \ref{pro:paatulos_ala} below are needed. Proposition \ref{pro:paatulos} gives the upper bound for the one-step error. It is an extension of Temam \cite{Temam} for the upper estimate and replaces the limit arguments by the inequality (\ref{eq:propo3ey}), which can be applied to any fixed time-net to get an upper bound for the approximation error. From Proposition \ref{pro:paatulos_ala} we get the lower bound for the one-step error. In the proof of Proposition \ref{pro:paatulos_ala} we use the same principal decomposition as in \cite{Temam}, but apply it to non-equidistant nets. We conclude the proof of our main result by considering multi-step error starting on page \pageref{page:paatulos_tod}.

\begin{prop}\label{pro:paatulos}
If for all $x \in E$
$$
\vits \frac{\partial^s}{\partial_{x_\beta}^q \partial_{x_\alpha}^r} \sigma_{ij}(x)\oits \leq C_3 \frac{Q_i(x)}{Q_\beta^q(x) Q_\alpha^r(x)}, \ q+r = s,\ q,r,s\in \vaalto 0,1,2\oaalto,
$$
$\vits b_i(x)\oits \leq C_3 Q_i(x)$ and $ A_{ii}(x) \geq \frac{1}{C_3}Q_i^2(x)$
for $i\in \vaalto 1,\ldots, d\oaalto$ and for some $C_3 > 0$,
then for $0\leq a \leq u < T$ it holds
\begin{align}\label{eq:propo3ey}
& \sum_{l=1}^d\sum_{k=1}^d \odotu  \vs \deri{x_k} F(u,X_u) - \deri{x_k} F(a,X_a)\os^2 \sigma_{k l}(X_u)^2\\
 &\leq D_3  \int_a^u \sup_{\alpha,\beta}
\odotu \vhaka A_{\alpha \alpha} (X_t) A_{\beta \beta}(X_t) \vits \deeri[2]{x_\alpha x_\beta} F(t,X_t)\oits^2 \ohaka dt,\nonumber
\end{align}
where $D_3 > 0$ depends at most on $C_3,d$ and $T$.
\end{prop}

\myproof
To keep the notation simple, we allow in the following that the constant $C>0$ may change from line to line.

Set
$$
v_a:= \vs\deri{x_k}F(a,X_a)\os_{k=1}^d \text{ and }
\phi_{kl}(u,x):=\vs\deri{x_k} F(u,x) - v_a^k\os\sigma_{kl}(x).
$$
Using this notation the assertion can be re-written as
$$
\sum_{l=1}^d \sum_{k=1}^d \odotu \phi_{kl}^2(u,X_u) \leq D  \int_a^u \sup_{\alpha,\beta}
\odotu \vhaka A_{\alpha \alpha} (X_t) A_{\beta \beta}(X_t) \vits \deeri[2]{x_\alpha x_\beta} F(t,X_t)\oits^2 \ohaka dt.
$$
By the definition of $\phi_{kl}$ we have that
\begin{align*}
\sum_{l=1}^d\sum_{k=1}^d \phi_{kl}^2(u,x) &=\sum_{l=1}^d \sum_{k=1}^d \vs \deri{x_k} F(u,x) - v_a^k\os^2 \sigma_{kl}^2(x)\\
&= \sum_{k=1}^d \vhaka \vs \deri{x_k} F(u,x) - v_a^k\os^2 \sum_{l=1}^d \sigma_{kl}^2(x) \ohaka.
\end{align*}
The assumptions on $\sigma$ give that 
$$
\frac{Q_k^2(x)}{C_3} \leq A_{kk}(x) = \sum_{l=1}^d \sigma_{kl}^2(x)  \leq d\ C_3 Q_k^2(x).
$$
This implies the equivalence
\begin{align}\label{eq:ekvivalenssi}
\frac{1}{C_3} \sum_{k=1}^d\vhaka \vs \deri{x_k} F(u,x) - v_a^k\os^2 Q_k^2(x)\ohaka &\leq
\sum_{l=1}^d\sum_{k=1}^d \phi_{kl}^2(u,x)\\
&\leq d\ C_3 \sum_{k=1}^d\vhaka \vs \deri{x_k} F(u,x) - v_a^k\os^2 Q_k^2(x)\ohaka. \nonumber
\end{align}
Lemma \ref{lemma:rajoittuneisuus} allows us to use the stopping argument from Lemma \ref{lemma:pysaytys}, which implies that
\begin{equation}\label{eq:odotus_phi_toiseen}
\odotu \phi_{kl}^2(u,X_u) = \int_a^u \odotu \vs \oA \phi_{kl}^2\os (v,X_v)dv + \sum_{m=1}^d \int_a^u \odotu\vs \deri{x_m} \phi_{kl}^2(v,X_v)\os b_m(X_v)dv.
\end{equation}
To prove our theorem we need to compute an upper bound for $\oA \phi_{kl}^2 (u,x)$ and for $\deri{x_m} \phi_{kl}^2(u,x) b_m(x)$. First we consider the term $\oA \phi_{kl}^2$:
\begin{align}\label{eq:A_phi_toiseen}
&\vits \oA \phi_{kl}^2(u,x)\oits \\\nonumber
&= \vits 2\phi_{kl}(u,x)\vs\oA \phi_{kl}\os(u,x) + \sum_{\alpha,\beta=1}^d \sum_{j=1}^d \vs\sigma_{\alpha j}(x) \deri{x_\alpha}\phi_{kl}(u,x)\os \vs \sigma_{\beta j}(x)\deri{x_\beta}\phi_{kl}(u,x)\os\oits \\\nonumber
%&\leq  \phi_{kl}^2(u,x) +  \vs \oA \phi_{kl}(u,x)\os^2 +\\\nonumber
%&\hspace*{2em} + \frac{1}{2} \sum_{\alpha,\beta = 1}^d \sum_{j=1}^d \vhaka \vs \sigma_{\alpha j}(x)\deri{x_\alpha}\phi_{kl}(u,x)\os^2 + \vs \sigma_{\beta j}(x)\deri{x_\beta}\phi_{kl}(u,x) \os^2 \ohaka\\\nonumber
%&=  \phi_{kl}^2(u,x) +  \vs \oA \phi_{kl}\os^2(u,x)
%+ d \sum_{\alpha=1}^d \sum_{j=1}^d \vs \sigma_{\alpha j}(x)\deri{x_\alpha}\phi_{kl}(u,x)\os^2. \nonumber
&\leq  \phi_{kl}^2(u,x) +  \vs \oA \phi_{kl}\os^2(u,x)
+ d \sum_{\alpha=1}^d \sum_{j=1}^d \vs \sigma_{\alpha j}(x)\deri{x_\alpha}\phi_{kl}(u,x)\os^2. \nonumber
\end{align}
Hence equation (\ref{eq:odotus_phi_toiseen}) implies that
\begin{align*}
\odotu \phi_{kl}^2(u,X_u) &\leq \int_a^u\odotu \phi_{kl}^2(v,X_v)dv + \int_a^u \odotu \vs \oA \phi_{kl}\os^2(v,X_v) dv\\
&+  d \sum_{\alpha=1}^d \sum_{j=1}^d \int_a^u \odotu \vs \sigma_{\alpha j}(X_v)\deri{x_\alpha}\phi_{kl}(v,X_v)\os^2 dv\\
&+ \sum_{m=1}^d \int_a^u \odotu \vits \vs\deri{x_m} \phi_{kl}^2(v,X_v)\os b_m(X_v) \oits dv,
\end{align*}
where the right-hand side is finite because of Lemma \ref{lemma:rajoittuneisuus}.
From Gronwall's Lemma (see Theorem \ref{the:gronwall} in Appendix) it follows that
\begin{align}\label{eq:equ1}
\odotu \phi_{kl}^2(u,X_u) &\leq \vhaka \int_a^u \odotu \vs \oA \phi_{kl}\os^2(v,X_v) dv +\right.\\\nonumber
&\hspace*{2em} + d \sum_{\alpha=1}^d \sum_{j=1}^d \int_a^u \odotu \vs \sigma_{\alpha j}(X_v)\deri{x_\alpha}\phi_{kl}(v,X_v)\os^2 dv\\\nonumber
&\hspace*{2em}+ \left. \sum_{m=1}^d \int_a^u \odotu \vits \deri{x_m} \phi_{kl}^2(v,X_v) b_m(X_v) \oits dv\ohaka
 e^{(u-a)}.\nonumber
\end{align}
To continue we need to find an upper bound for the above expression. We start with $\oA \phi_{kl}$ and have, by definition, that
\begin{align*}
\oA \phi_{kl}(u,x)
&=\vs\deeri[2]{t x_k} F(u,x)\os \sigma_{kl}(x)+\\
&\hspace*{2em} +\frac{1}{2}  \sum_{\alpha,\beta=1}^d A_{\alpha \beta}(x)\vhaka \vs \deeri[2]{x_\alpha x_\beta} \sigma_{kl}(x)\os \vs \deri{x_k} F(u,x) - v_a^k\os \right.\\
& \hspace*{10em}+\vs \deri{x_\beta}\sigma_{kl}(x)\os \vs \deeri[2]{x_\alpha x_k} F(u,x)\os\\
& \hspace*{10em}+\vs \deri{x_\alpha}\sigma_{kl}(x)\os \vs \deeri[2]{x_\beta x_k} F(u,x)\os\\
& \hspace*{10em}+ \left. \sigma_{kl}(x) \vs \deeri[3]{x_\alpha x_\beta x_k} F(u,x)\os\ohaka.
\end{align*}
Taking the derivative with respect to $x_k$ in the partial differential equation (\ref{eq:bpde}) we get that
\begin{align*}
\deri{x_k} &\deri{t}F(u,x) +\frac{1}{2} \sum_{\alpha, \beta=1}^d A_{\alpha\beta}(x) \deeri[3]{x_k x_\alpha x_\beta}F(u,x)\\
&= - \frac{1}{2} \sum_{\alpha, \beta=1}^d\vs \deri{x_k} A_{\alpha \beta}(x)\os  \deeri[2]{x_\alpha x_\beta}F(u,x).
\end{align*}
Now we can replace the derivative with respect to $t$ and the third order derivatives in the formula for $\oA \phi_{kl}(u,x)$ by second order derivatives:
\begin{align*}
\oA \phi_{kl}(u,x) &=
\frac{1}{2}  \sum_{\alpha,\beta=1}^d \vhaka A_{\alpha \beta}(x)\vs \deeri[2]{x_\alpha x_\beta} \sigma_{kl}(x)\os \vs \deri{x_k} F(u,x) - v_a^k\os \right.\\
& \hspace*{8em}+ A_{\alpha \beta}(x) \vs \deri{x_\beta}\sigma_{kl}(x)\os \vs \deeri[2]{x_\alpha x_k} F(u,x)\os\\
& \hspace*{8em}+ A_{\alpha \beta}(x) \vs \deri{x_\alpha}\sigma_{kl}(x)\os \vs \deeri[2]{x_\beta x_k} F(u,x)\os\\
& \hspace*{8em}- \left. \sigma_{kl}(x) \vs \deri{x_k} A_{\alpha \beta}(x)\os \vs \deeri[2]{x_\alpha x_\beta} F(u,x)\os\ohaka.
\end{align*}
It follows from the definition of the matrix $A$ and the assumption on the matrix $\sigma$ that
$$
\vits A_{\alpha \beta}(x)\oits \leq C Q_\alpha(x) Q_\beta(x)
$$
and 
$$
\vits \deri{x_k}A_{\alpha \beta} (x)\oits \leq C \frac{Q_\alpha(x) Q_\beta(x)}{Q_k(x)}.
$$
Now we can bound the function $\vs \oA\phi_{kl}\os^2(u,x)$ from the above by
\begin{align}\label{eq:equ2}
\vs \oA\phi_{kl}\os^2(u,x) &\leq
C\sum_{\alpha,\beta=1}^d\vhaka \vs \deri{x_k} F(u,x) - v_a^k\os^2   Q_k^2(x)  \right.\\\nonumber
& \hspace*{5em}+  Q_\alpha^2(x) Q_k^2(x) \vs \deeri[2]{x_\alpha x_k} F(u,x)\os^2\\\nonumber
& \hspace*{5em}+  Q_\beta^2(x) Q_k^2(x) \vs \deeri[2]{x_\beta x_k} F(u,x)\os^2\\\nonumber
& \hspace*{5em}+ \left. Q_\alpha^2(x) Q_\beta^2(x) \vs \deeri[2]{x_\alpha x_\beta} F(u,x)\os^2\ohaka.\nonumber
\end{align}
%\vfill\eject
For $ \vs \sigma_{\alpha j}(x) \deri{x_\alpha} \phi_{kl}(u,x)\os^2$ we get that
\begin{align}\label{eq:equ3}
\nonumber
&\vs \sigma_{\alpha j}(x) \deri{x_\alpha} \phi_{kl}(u,x)\os^2\\
&\leq C \vs Q_k^2(x)\vs \deri{x_k}F(u,x) - v_a^k\os^2 + Q_{\alpha}^2(x) Q_{k}^2(x) \vs\deeri[2]{x_\alpha x_k}F(u,x)\os^2\os.\\\nonumber
\end{align}
The term including $b_m$ can be bounded as follows:
\begin{align}\label{eq:equ4}
\nonumber
&\vits\vs\deri{x_m} \phi_{kl}^2(u,x)\os b_m(x)\oits\\\nonumber
&\leq 2 \vits \deri{x_k} F(u,x) - v_a^k\oits \vits\sigma_{kl}(x)\oits \vits b_m(x)\oits \times\\\nonumber
&\hspace*{5em} \times  \vs \vits \deri{x_m} \sigma_{kl}(x)\vs \deri{x_k}F(u,x) - v_a^k\os\oits + \vits \sigma_{kl}(x) \oits \vits \deeri[2]{x_m x_k}F(u,x)\oits\os \\\nonumber
&\leq C \vits \deri{x_k} F(u,x) - v_a^k\oits Q_k(x) Q_m(x)\times\\\nonumber
&\hspace*{2em} \times  \vs \frac{Q_k(x)}{Q_m(x)} \vits \deri{x_k}F(u,x) - v_a^k\oits + Q_k(x) \vits \deeri[2]{x_m x_k}F(u,x)\oits\os \\\nonumber
&\leq C  \vs \vits \deri{x_k} F(u,x) - v_a^k\oits^2 Q_k^2(x) + Q_k^2(x)Q_m^2(x) \vits \deeri[2]{x_m x_k}F(u,x)\oits^2\os,\\
\end{align}
where we used that
\begin{align*}
&\vits \deri{x_k} F(u,x) - v_a^k\oits  Q_k^2(x)Q_m(x) \vits \deeri[2]{x_m x_k} F(u,x)\oits\\
& \leq \vits \deri{x_k} F(u,x) - v_a^k\oits^2 Q_k^2(x) + Q_k^2(x)Q_m^2(x) \vits \deeri[2]{x_m x_k} F(u,x)\oits^2.
\end{align*}
Now the expectation of $\phi_{kl}^2(u,X_u)$ can be bounded by
\begin{align*}
\odotu \phi_{kl}^2(u,X_u)
&\leq C\int_a^u  \odotu\vs \deri{x_k} F(v,X_v) - v_a^k\os^2 Q_k^2(X_v)dv\\
& + C \int_a^u \sup_{\alpha,\beta}\odotu Q_\alpha^2(X_v) Q_\beta^2(X_v) \vs \deeri[2]{x_\alpha x_\beta} F(v,X_v)\os^2dv,
\end{align*}
where we use (\ref{eq:equ1}), (\ref{eq:equ2}), (\ref{eq:equ3}) and (\ref{eq:equ4}).
From the above and (\ref{eq:ekvivalenssi}) we get
\begin{align*}
&\sum_{k=1}^d \odotu\vhaka \vs \deri{x_k} F(u,X_u) - v_a^k\os^2 Q_k^2(X_u)\ohaka \leq C \sum_{l=1}^d\sum_{k=1}^d \odotu \phi_{kl}^2(u,X_u)\\
&\leq C  \int_a^u \odotu \sum_{k=1}^d  \vs \deri{x_k} F(v,X_v) - v_a^k\os^2 Q_k^2(X_v) dv +\\
&\hspace*{2em} + C \int_a^u \sup_{\alpha,\beta} \odotu Q_\alpha^2(X_v) Q_\beta^2(X_v) \vs \deeri[2]{x_\alpha x_\beta} F(v,X_v)\os^2 dv.
\end{align*}
Gronwall's lemma (Theorem \ref{the:gronwall}) gives
\begin{align}\label{eq:paatulos_tod_loppu}
&\sum_{k=1}^d \odotu\vhaka \vs \deri{x_k} F(u,X_u) - v_a^k\os^2 Q_k^2(X_u)\ohaka\\\nonumber
&\leq  e^{C(u-a)} C \int_a^u \sup_{\alpha,\beta} \odotu Q_\alpha^2(X_v) Q_\beta^2(X_v) \vs \deeri[2]{x_\alpha x_\beta} F(v,X_v)\os^2 dv\nonumber
\end{align}
and the assertion follows from (\ref{eq:ekvivalenssi}).
\stopproof

\begin{prop}\label{pro:paatulos_ala}
If for all $x \in E$
$$
\vits \frac{\partial^s}{\partial_{x_\beta}^q \partial_{x_\alpha}^r} \sigma_{ij}(x)\oits \leq C_4 \frac{Q_i(x)}{Q_\beta^q(x) Q_\alpha^r(x)}, \ q+r = s,\ q,r,s\in \vaalto 0,1,2\oaalto,
$$
and
$$
 \vits b_i(x)\oits \leq C_4 Q_i(x)
$$
for $i\in \vaalto 1,\ldots, d\oaalto$ and for some $C_4 > 0$, then for $0\leq a < t < T$ it holds that
\begin{align*}
&\sum_{j=1}^d \odotu  \vs \sum_{i=1}^d\vs \deri{x_i}F(t,X_t) - \deri{x_i}F(a,X_a)\os \sigma_{ij}(X_u) \os^2 
\geq  \int_a^t H^2(u) du\\
&  - D_4\int_a^t \vhaka \vs \int_a^u \sup_{\alpha,\beta}\odotu Q_\alpha^2 (X_v)
Q_\beta^2(X_v) \vits \deeri[2]{x_\alpha
x_\beta}F(v,X_v)\oits^2 dv\os^\frac{1}{2}\times \right.\\
&\hspace*{5em} \times \vs \sup_{\alpha,\beta} \odotu Q_\alpha^2 (X_u)
Q_\beta^2(X_u) \vits \deeri[2]{x_\alpha
x_\beta}F(u,X_u)\oits^2\os^\frac{1}{2} +\\
&\hspace*{5em} + \left. \int_a^u \sup_{\alpha,\beta}\odotu Q_\alpha^2 (X_v)
Q_\beta^2(X_v) \vits \deeri[2]{x_\alpha
x_\beta}F(v,X_v)\oits^2 dv \ohaka du,
\end{align*}
where the function $H$ is defined in Theorem \ref{the:paatulos} and $D_4 > 0$ depends at most on $C_4,d$ and $T$.
\end{prop}

\myproof

To abbreviate the notation we assume again that $C>0$ may change from line to line. We let
$$
\phi_u^{ij}:=\vs \deri{x_i}F(u,X_u) - v_a^i\os \sigma_{ij}(X_u),
$$
where $u\in [a,t]$ and $v_a^i:= \deri{x_i}F(a,X_a)$.
It\^o's formula gives that, a.s.,
\begin{align*}
\phi_t^{ij} &= \int_a^t \deeri[2]{t x_i}F(u,X_u) \sigma_{ij}(X_u) du\\
&\hspace*{1em} + \sum_{\alpha=1}^d \int_a^t b_\alpha(X_u) \vhaka \deeri[2]{x_\alpha x_i}F(u,X_u) \sigma_{ij}(X_u) \right.\\
&\hspace*{7em} \left. +\vs \deri{x_i}F(u,X_u) - v_a^i \os \deri{x_\alpha}\sigma_{ij}(X_u) \ohaka du\\
&\hspace*{1em} + \sum_{\alpha=1}^d \sum_{n=1}^d \int_a^t \vhaka \deeri[2]{x_\alpha x_i}F(u,X_u) \sigma_{ij}(X_u) \right.\\
&\hspace*{7em} \left. +\vs \deri{x_i}F(u,X_u) - v_a^i \os \deri{x_\alpha}\sigma_{ij}(X_u) \ohaka\sigma_{\alpha n}(X_u)dW_u^n\\
&\hspace*{1em} + \frac{1}{2} \sum_{\alpha, \beta=1}^d \int_a^t\vhaka \deeri[3]{x_\beta x_\alpha x_i}F(u,X_u) \sigma_{ij}(X_u) + \deeri[2]{x_\alpha x_i} F(u,X_u) \deri{x_\beta} \sigma_{ij}(X_u) +\right.\\
&\hspace*{6em}+ \deeri[2]{x_\beta x_i}F(u,X_u) \deri{x_\alpha}\sigma_{ij}(X_u) +\\
&\hspace*{6em}+ \left. \vs \deri{x_i}F(u,X_u) - v_a^i\os \deeri[2]{x_\beta x_\alpha} \sigma_{ij}(X_u) \ohaka A_{\alpha \beta}(X_u) du.
\end{align*}
From the above we deduce
\begin{align*}
&d\langle \phi^{ij},\phi^{kj}\rangle_u\\
&= \sum_{\alpha=1}^d \sum_{\beta=1}^d \vhaka \vs \deeri[2]{x_\alpha x_i} F(u,X_u) \sigma_{ij}(X_u) + \vs \deri{x_i} F(u,X_u) - v_a^i\os \deri{x_\alpha} \sigma_{ij}(X_u) \os \times \right.\\
&\times \left. \vs \deeri[2]{x_\beta x_k}F(u,X_u) \sigma_{kj}(u,X_u) + \vs \deri{x_k} F(u,X_u) - v_a^k\os \deri{x_\beta} \sigma_{kj}(X_u)\os A_{\alpha \beta}(X_u) \ohaka du
\end{align*}
and, using the equality $\oA F = 0$ (cf. (\ref{eq:oA}) and (\ref{eq:bpde})), we get that
\begin{align*}
&d\phi_u^{ij}\\
&= \sum_{\alpha=1}^d \sum_{n=1}^d \vs \deeri[2]{x_\alpha x_i}F(u,X_u) \sigma_{ij}(X_u) + \vs \deri{x_i}F(u,X_u) - v_a^i\os \deri{x_\alpha}\sigma_{ij}(X_u)\os \sigma_{\alpha n}(X_u) dW_u^n\\
&+ \frac{1}{2} \sum_{\alpha,\beta=1}^d A_{\alpha \beta}(X_u) \vhaka \deeri[2]{x_\alpha x_i} F(u,X_u) \deri{x_\beta}\sigma_{ij}(X_u) + \deeri[2]{x_\beta x_i} F(u,X_u) \deri{x_\alpha}\sigma_{ij}(X_u)\right.\\
&\hspace*{6em} +\left. \vs \deri{x_i}F(u,X_u) - v_a^i\os \deeri[2]{x_\beta x_\alpha} \sigma_{ij}(X_u)\ohaka du\\
&- \frac{1}{2}\sum_{\alpha,\beta=1}^d \vhaka\deri{x_i} A_{\alpha \beta}(X_u) \deeri[2]{x_\beta x_\alpha}F(u,X_u) \sigma_{ij}(X_u)\ohaka du\\
&+ \sum_{\alpha=1}^d b_\alpha(X_u) \vhaka \deeri[2]{x_\alpha x_i}F(u,X_u) \sigma_{ij}(X_u) \right.\\
&\hspace*{7em} \left. +\vs \deri{x_i}F(u,X_u) - v_a^i \os \deri{x_\alpha}\sigma_{ij}(X_u) \ohaka du.
\end{align*}
Let
\begin{equation*}
S_N := \inf\vaalto u\geq a: \vits \phi_u^{ij} \oits \geq N \text{ or } \vnorm X_u\onorm \geq N \text{ for some } i,j \in \vaalto 1,\ldots, d\oaalto \oaalto \wedge t
\end{equation*}
for $N = 1,2,\ldots$. Because of $S_N \nearrow t$ a.s. as $N\to \infty$ and Lemma \ref{lemma:rajoittuneisuus} one has that
\begin{equation*}
\sum_{j=1}^d \odotu \vs \sum_{i=1}^d \phi_t^{ij}\os^2 = \lim_{N\to\infty} \sum_{j=1}^d \odotu \vs \sum_{i=1}^d \phi_{S_N}^{ij}\os^2.
\end{equation*}
Using the integration by parts formula for semi-martingales, we get 
$$
\sum_{j=1}^d \odotu \vs \sum_{i=1}^d \phi_{S_N}^{ij}\os^2 = \sum_{j=1}^d 2 \odotu \int_a^{S_N} \sum_{i,k=1}^d \phi_u^{kj} d\phi_u^{ij} + \sum_{i,k=1}^d \odotu \int_a^{S_N}  d\langle \phi^{ij},\phi^{kj}\rangle_u.
$$
Because of the choice of the stopping time $S_N$, the expected value of the "$dW_u^n$-terms" vanishes. For the rest, we obtain as main term
$$
\odotu  \int_a^{S_N} \sum_{\alpha,\beta,i,k=1}^d
A_{\alpha \beta}(X_u) A_{i k}(X_u) \deeri[2]{x_\alpha x_i} F(u,X_u)  \deeri[2]{x_\beta x_k} F(u,X_u) du 
$$
and (after some computation) terms of the type
$$
\odotu  \int_a^{S_N}  A_{\alpha \beta}(X_u) \sigma_{kj}(X_u) \deri{x_\beta}\sigma_{ij}(X_u) \vs \deri{x_k}F(u,X_u) - v_a^k \os
\deeri[2]{x_\alpha x_i}F(u,X_u) du
$$
and
$$
\odotu  \int_a^{S_N} \vs \deri{x_k} F(u,X_u) - v_a^k \os \sigma_{kj}(X_u) b_\alpha(X_u) \vs \deri{x_i}F(u,X_u) - v_a^i \os \deri{x_\alpha}\sigma_{ij}(X_u)  du.
$$
Using the assumptions on the matrix $\sigma$ and the vector $b$ we can bound these terms by
\begin{align*}
&\vits  A_{\alpha \beta}(X_u) \sigma_{kj}(X_u) \deri{x_\beta}\sigma_{ij}(X_u) \vs \deri{x_k}F(u,X_u) - v_a^k \os
\deeri[2]{x_\alpha x_i}F(u,X_u) \oits\\
%&\leq C Q_\alpha(X_u) Q_\beta(X_u) Q_k(X_u) \frac{Q_i(X_u)}{Q_\beta(X_u)} \vits \deri{x_k}F (u,X_u)- v_a^k \oits \vits \deeri[2]{x_\alpha x_i} F(u,X_u) \oits\\
&\leq C   Q_k(X_u)  \vits \deri{x_k}F(u,X_u) - v_a^k \oits Q_\alpha(X_u) Q_i(X_u) \vits \deeri[2]{x_\alpha x_i} F(u,X_u) \oits
\end{align*}
and 
\begin{align*}
&\vits \vs \deri{x_k} F(u,X_u) - v_a^k \os \sigma_{kj}(X_u) b_\alpha(X_u) \vs \deri{x_i}F(u,X_u) - v_a^i \os \deri{x_\alpha}\sigma_{ij}(X_u) \oits\\
&\leq C Q_k(X_u) \vits \deri{x_k}F(u,X_u) - v_a^k \oits Q_i(X_u) \vits \deri{x_i}F(u,X_u) - v_a^i
\oits.
\end{align*}
H\"older's inequality and (\ref{eq:paatulos_tod_loppu}) give
\begin{align*}
&\odotu Q_k(X_u)  \vits \deri{x_k}F(u,X_u) - v_a^k \oits Q_\alpha(X_u) Q_i(X_u) \vits \deeri[2]{x_\alpha x_i} F(u,X_u)
\oits\\
&\leq \vs \odotu Q_k^2(X_u) \vits \deri{x_k}F(u,X_u) -
v_a^k\oits^2 \os^\frac{1}{2} \vs \sup_{\alpha', \beta'}\odotu Q_{\alpha'}^2 (X_u)
Q_{\beta'}^2(X_u) \vits \deeri[2]{x_{\alpha'}
x_{\beta'}}F(u,X_u)\oits^2\os^\frac{1}{2}\\
&\leq C\vs \int_a^u \sup_{\alpha',\beta'}\odotu Q_{\alpha'}^2 (X_v)
Q_{\beta'}^2(X_v) \vits \deeri[2]{x_{\alpha'}
x_{\beta'}}F(v,X_v)\oits^2 dv\os^\frac{1}{2}\times\\
&\hspace*{2em}\times \vs\sup_{\alpha',\beta'} \odotu Q_{\alpha'}^2 (X_u)
Q_{\beta'}^2(X_u) \vits \deeri[2]{x_{\alpha'}
x_{\beta'}}F(u,X_u)\oits^2\os^\frac{1}{2}.
\end{align*}
Moreover
\begin{align*}
&\odotu Q_k(X_u)  \vits \deri{x_k}F(u,X_u) - v_a^k \oits Q_i(X_u)  \vits \deri{x_i}F(u,X_u) - v_a^i \oits\\
&\leq \vs \odotu Q_k^2(X_u) \vits \deri{x_k}F(u,X_u) -
v_a^k\oits^2 \os^\frac{1}{2} \vs \odotu Q_i^2(X_u) \vits \deri{x_i}F(u,X_u) -
v_a^i\oits^2 \os^\frac{1}{2}\\
&\leq C  \int_a^u \sup_{\alpha,\beta}\odotu Q_\alpha^2 (X_v)
Q_\beta^2(X_v) \vits \deeri[2]{x_\alpha
x_\beta}F(v,X_v)\oits^2 dv
\end{align*}
and the assertion follows by $N\to \infty$.
\stopproof

{\bf Proof of Theorem \ref{the:paatulos}.}\\
\label{page:paatulos_tod}
Also in this proof, we use the same notation for different constants.
First we consider the upper bound for the approximation error. Let $\eps \in (0,T)$. Using Doob's inequality together with H\"older's inequality we see that
\begin{align*}
&\vs \odotu \sup_{t\in[0,T-\eps]} \vits \sum_{i=1}^n \sum_{k=1}^d \int_{t_{i-1}\wedge t}^{t_i\wedge t} \vs\deri{x_k} F(u,X_u) - \deri{x_k}F(t_{i-1},X_{t_{i-1}})\os dX_u^k \oits^2\os^\frac{1}{2}\\
&\leq \vs \odotu \sup_{t\in[0,T-\eps]} \vits \sum_{i=1}^n \sum_{k=1}^d \int_{t_{i-1}\wedge t}^{t_i\wedge t} \vs \deri{x_k} F(u,X_u) - \deri{x_k}F(t_{i-1},X_{t_{i-1}}) \os b_k(X_u)du\oits^2\os^\frac{1}{2} + \\
&+ \vs \odotu \sup_{t\in[0,T-\eps]} \vits \sum_{i=1}^n \sum_{l=1}^d \int_{t_{i-1}\wedge t}^{t_i\wedge t}  \sum_{k=1}^d\vs \deri{x_k} F(u,X_u) - \deri{x_k}F(t_{i-1},X_{t_{i-1}})\os \sigma_{kl}(X_u) dW_u^l\oits^2 \os^\frac{1}{2} \\
&\leq  \vs \odotu \vits \sum_{i=1}^n \sum_{k=1}^d \int_{t_{i-1} \wedge T-\eps}^{t_i \wedge T-\eps} \vits \deri{x_k} F(u,X_u) - \deri{x_k}F(t_{i-1},X_{t_{i-1}}) \oits \vits b_k(X_u)\oits du\oits^2\os^\frac{1}{2} + \\
&+ 2 \vs \odotu \vits \sum_{i=1}^n \sum_{l=1}^d \int_{t_{i-1} \wedge T-\eps}^{t_i \wedge T-\eps}  \sum_{k=1}^d\vs \deri{x_k} F(u,X_u) - \deri{x_k}F(t_{i-1},X_{t_{i-1}})\os \sigma_{kl}(X_u) dW_u^l\oits^2 \os^\frac{1}{2} \\
&\leq \sum_{k=1}^d \sqrt{T} \vs \odotu \sum_{i=1}^n \int_{t_{i-1}}^{t_i} \vits \deri{x_k} F(u,X_u) - \deri{x_k}F(t_{i-1},X_{t_{i-1}}) \oits^2 \vits b_k(X_u)\oits^2 du \os^\frac{1}{2} + \\
&+ 2 \vs \odotu \vits \sum_{i=1}^n \sum_{l=1}^d \int_{t_{i-1} \wedge T-\eps}^{t_i \wedge T-\eps}  \sum_{k=1}^d\vs \deri{x_k} F(u,X_u) - \deri{x_k}F(t_{i-1},X_{t_{i-1}})\os \sigma_{kl}(X_u) dW_u^l\oits^2 \os^\frac{1}{2} \\
&:= B_1 + B_2.
\end{align*}
Inequality (\ref{eq:paatulos_tod_loppu}) gives that 
\begin{align*}
B_1 &\leq C_1 \sqrt{T} \sum_{k=1}^d \vs\sum_{i=1}^n \int_{t_{i-1}}^{t_i} \odotu \vs \deri{x_k} F(u,X_u) - \deri{x_k}F(t_{i-1},X_{t_{i-1}}) \os^2 Q_k^2(X_u)du\os^\frac{1}{2}\\
&\leq \sqrt{T} C  \sum_{k=1}^d  \vs \sum_{i=1}^n \int_{t_{i-1}}^{t_i} \int_{t_{i-1}}^u \sup_{\alpha, \beta} \odotu Q^2_\alpha (X_v) Q^2_\beta(X_v) \vs \deeri[2]{x_\alpha x_\beta} F(v,X_v)\os^2 dv du\os^\frac{1}{2}.
\end{align*}
For $B_2$ the It\^o-isometry and the orthogonality of stochastic integrals give
\begin{align*}
B_2^2 &= 4 \odotu \vits \sum_{i=1}^n \sum_{l=1}^d \int_{t_{i-1} \wedge T-\eps}^{t_i \wedge T-\eps}  \sum_{k=1}^d\vs \deri{x_k} F(u,X_u) - \deri{x_k}F(t_{i-1},X_{t_{i-1}})\os \sigma_{kl}(X_u) dW_u^l\oits^2  \\
&= 4 \sum_{i=1}^n \sum_{l=1}^d \odotu \vits \int_{t_{i-1} \wedge T-\eps}^{t_i \wedge T-\eps}  \sum_{k=1}^d\vhaka\deri{x_k} F(u,X_u) - \deri{x_k}F(t_{i-1},X_{t_{i-1}})\ohaka\sigma_{kl}(X_u) dW_u^l\oits^2 \\
&= 4 \sum_{i=1}^n \sum_{l=1}^d \int_{t_{i-1} \wedge T-\eps}^{t_i \wedge T-\eps} \odotu \vits \sum_{k=1}^d\vhaka\deri{x_k} F(u,X_u) - \deri{x_k}F(t_{i-1},X_{t_{i-1}})\ohaka\sigma_{kl}(X_u)\oits^2 du \\
&\leq 4 d \sum_{i=1}^n \sum_{l=1}^d \sum_{k=1}^d \int_{t_{i-1}}^{t_i}  \odotu  \vs \deri{x_k} F(u,X_u) - \deri{x_k}F(t_{i-1},X_{t_{i-1}})\os^2 \sigma_{k l}(X_u)^2 du.
\end{align*}
Letting $\eps \searrow 0$ we get by monotone convergence that
\begin{align*}
&\vs \odotu \sup_{t\in [0,T]} \vits \sum_{i=1}^n \sum_{k=1}^d \int_{t_{i-1}\wedge t}^{t_i\wedge t} \vs\deri{x_k} F(u,X_u) - \deri{x_k}F(t_{i-1},X_{t_{i-1}})\os dX_u^k \oits^2\os^\frac{1}{2}\\
& \leq 
\sqrt{T} C  \sum_{k=1}^d  \vs \sum_{i=1}^n \int_{t_{i-1}}^{t_i} \int_{t_{i-1}}^u \sup_{\alpha, \beta} \odotu Q^2_\alpha (X_v) Q^2_\beta(X_v) \vs \deeri[2]{x_\alpha x_\beta} F(v,X_v)\os^2 dv\os^\frac{1}{2} + \\
&\hspace*{2em}+ \vs4 d \sum_{i=1}^n \sum_{l=1}^d \sum_{k=1}^d \int_{t_{i-1}}^{t_i}  \odotu  \vs \deri{x_k} F(u,X_u) - \deri{x_k}F(t_{i-1},X_{t_{i-1}})\os^2 \sigma_{k l}(X_u)^2 du\os^\frac{1}{2}.
\end{align*}
The assertion for the upper bound follows from Proposition \ref{pro:paatulos} and Lemma \ref{tthe:geiss} in Appendix.

Now we continue with the lower bound of the approximation error. Let $[A,B]$ be a subinterval of $(r,s)$ such that
\begin{equation}\label{eq:vakio_C_B}
0< \frac{(B-A) C_B}{(T-B)^{2\theta}} \leq \frac{C_H}{4},
\end{equation}
where $C_H$ is taken from (\ref{eq:oletus_H}) and the constant $C_B > 0$ satisfies (cf. (\ref{eq:tulos2oletus}) and (\ref{eq:paatulos_tod_loppu}))
\begin{equation}\label{eq:drifti_raja}
d \sum_{k=1}^d\odotu \vs \deri{x_k}F(u,X_u) - \deri{x_k}F(a,X_a)\os^2 b_k^2(X_u) \leq C_B \int_a^u \frac{1}{(T-u)^{2\theta}}dv
\end{equation}
for $A \leq a < u \leq B$. 
Let us now consider the approximation error inside the interval $[A, B]$. Denote $I_n:= \vaalto i: A\leq t_{i-1}^n \leq t_i^n \leq B\oaalto$ and denote in both cases for the lower estimate (cf. Theorem \ref{the:paatulos} cases $(L1)$ and $(L2)$) the sequence of time-nets by $(t_i^n)_{i=0}^n$. Note that for large $n$ the set $I_n$ is not an empty set because, in both cases, we have that $\sup_{i=1,\ldots, n} (t_i^n - t_{i-1}^n) \leq C_\tau/n$. Now on $[A,B]$ we get that
\begin{align*}
&\vs \odotu \vits \sum_{i\in I_n} \sum_{k=1}^d \int_{t_{i-1}^n}^{t_i^n} \vs\deri{x_k} F(u,X_u) - \deri{x_k}F(t_{i-1}^n,X_{t_{i-1}^n})\os dX_u^k \oits^2\os^\frac{1}{2}\\
&\geq  \vits \vs \odotu  \vits \sum_{i\in I_n} \sum_{l=1}^d \int_{t_{i-1}^n}^{t_i^n}  \sum_{k=1}^d\vs \deri{x_k} F(u,X_u) - \deri{x_k}F(t_{i-1}^n,X_{t_{i-1}^n})\os \sigma_{kl}(X_u) dW_u^l\oits^2 \os^\frac{1}{2} - \right. \\
&- \left. \vs \odotu  \vits \sum_{i\in I_n} \sum_{k=1}^d \int_{t_{i-1}^n}^{t_i^n} \vs \deri{x_k} F(u,X_u) - \deri{x_k}F(t_{i-1}^n,X_{t_{i-1}^n}) \os b_k(X_u)du\oits^2\os^\frac{1}{2}\oits.
\end{align*}
This implies that
\begin{align*}
& \odotu \vits \sum_{i\in I_n} \sum_{k=1}^d \int_{t_{i-1}^n}^{t_i^n} \vs\deri{x_k} F(u,X_u) - \deri{x_k}F(t_{i-1}^n,X_{t_{i-1}^n})\os dX_u^k \oits^2\\
&\geq \frac{1}{2} \odotu  \vits \sum_{i\in I_n} \sum_{l=1}^d \int_{t_{i-1}^n}^{t_i^n}  \sum_{k=1}^d\vs \deri{x_k} F(u,X_u) - \deri{x_k}F(t_{i-1}^n,X_{t_{i-1}^n})\os \sigma_{kl}(X_u) dW_u^l\oits^2 - \\
&- \odotu \vits \sum_{i\in I_n} \sum_{k=1}^d  \int_{t_{i-1}^n}^{t_i^n} \vs \deri{x_k} F(u,X_u) - \deri{x_k}F(t_{i-1}^n,X_{t_{i-1}^n}) \os b_k (X_u)du\oits^2\\
&\geq \frac{1}{2} \odotu  \vits \sum_{i\in I_n} \sum_{l=1}^d \int_{t_{i-1}^n}^{t_i^n}  \sum_{k=1}^d\vs \deri{x_k} F(u,X_u) - \deri{x_k}F(t_{i-1}^n,X_{t_{i-1}^n})\os \sigma_{kl}(X_u) dW_u^l\oits^2\\
&- (B-A)   \sum_{i\in I_n} \ d \sum_{k=1}^d\odotu   \int_{t_{i-1}^n}^{t_i^n} \vs \deri{x_k} F(u,X_u) - \deri{x_k}F(t_{i-1}^n,X_{t_{i-1}^n}) \os^2 b_k^2(X_u)du.
\end{align*}
Now (\ref{eq:vakio_C_B}) and (\ref{eq:drifti_raja}) give that
\begin{align*}
&(B-A)   \sum_{i\in I_n} \ d \sum_{k=1}^d\odotu   \int_{t_{i-1}^n}^{t_i^n} \vs \deri{x_k} F(u,X_u) - \deri{x_k}F(t_{i-1}^n,X_{t_{i-1}^n}) \os^2 b_k^2(X_u)du\\
&\leq (B-A)   \sum_{i\in I_n} \int_{t_{i-1}^n}^{t_i^n} \int_{t_{i-1}^n}^{u} \frac{C_B}{(T-v)^{2\theta}} dvdu\\
&\leq \sum_{i\in I_n} \int_{t_{i-1}^n}^{t_i^n} \int_{t_{i-1}^n}^{u} \frac{H^2(v)}{4}dv du.
\end{align*}
Let us now consider the lower bound for 
$$
\vs \odotu  \vits \sum_{i\in I_n} \sum_{l=1}^d \int_{t_{i-1}^n}^{t_i^n}  \sum_{k=1}^d\vs \deri{x_k} F(u,X_u) - \deri{x_k}F(t_{i-1}^n,X_{t_{i-1}^n})\os \sigma_{kl}(X_u) dW_u^l\oits^2 \os^\frac{1}{2}.
$$
Using the It\^o-isometry, for $0\leq a < b < T$, we get that
\begin{align*}
B_3:=&\odotu\vits \sum_{l=1}^d \int_a^b \sum_{k=1}^d \vs \deri{x_k}F(t,X_t) - \deri{x_k}F(a,X_a)\os \sigma_{kl}(X_t) dW_t^l \oits^2 \\
&= \sum_{l=1}^d \odotu \int_a^b \vs \sum_{k=1}^d\vs \deri{x_k}F(t,X_t) - \deri{x_k}F(a,X_a)\os \sigma_{kl}(X_t) \os^2 dt.
\end{align*}
Assuming $b-a \leq 1$, Proposition \ref{pro:paatulos_ala} together with (\ref{eq:tulos2oletus}) imply that
\begin{align*}
B_3 &\geq  \int_a^b \int_a^t H^2(u) du dt
 - D_4 \int_a^b \int_a^t \vs \int_a^u
\frac{1}{(T-v)^{2\theta}}dv\os^\frac{1}{2} \vs
\frac{1}{(T-u)^{2\theta}}\os^\frac{1}{2} du dt\\
& - D_4 \int_a^b \int_a^t \int_a^u
\frac{1}{(T-v)^{2\theta}}dv du dt\\
&\geq \int_a^b \int_a^t H^2(u) du dt
 - 2 D_4 \sqrt{b-a} \int_a^b \int_a^t \frac{1}{(T-u)^{2\theta}} du dt.
\end{align*}
Considering the multi-step error for the approximation we get that
\begin{align*}
& \odotu \vits \sum_{i\in I_n} \sum_{k=1}^d \int_{t_{i-1}^n}^{t_i^n} \vs\deri{x_k} F(u,X_u) - \deri{x_k}F(t_{i-1}^n,X_{t_{i-1}^n})\os dX_u^k \oits^2\\
&\geq \frac{1}{2} \odotu  \vits \sum_{i\in I_n} \sum_{l=1}^d \int_{t_{i-1}^n}^{t_i^n}  \sum_{k=1}^d\vs \deri{x_k} F(u,X_u) - \deri{x_k}F(t_{i-1}^n,X_{t_{i-1}^n})\os \sigma_{kl}(X_u) dW_u^l\oits^2\\
& - (B-A)  \sum_{i\in I_n} \ d \sum_{k=1}^d \odotu   \int_{t_{i-1}^n}^{t_i^n} \vs \deri{x_k} F(u,X_u) - \deri{x_k}F(t_{i-1}^n,X_{t_{i-1}^n}) \os^2 b_k^2(X_u)du \\
&\geq \frac{1}{2} \sum_{i\in I_n} \vhaka \int_{t_{i-1}^n}^{t_i^n} \int_{t_{i-1}^n}^t H^2(u)dudt - C \int_{t_{i-1}^n}^{t_i^n} \int_{t_{i-1}^n}^t  \frac{\sqrt{t_i^n- t_{i-1}^n}}{(T-u)^{2\theta}} dudt\ohaka\\
& - \sum_{i\in I_n} \int_{t_{i-1}^n}^{t_i^n} \int_{t_{i-1}^n}^t \frac{H^2(u)}{4}dudt\\
&\geq \sum_{i\in I_n}  \int_{t_{i-1}^n}^{t_i^n} (t_i^n - t) \frac{H^2(t)}{4}dt - \frac{C}{\sqrt{n}} \sum_{i\in I_n} \int_{t_{i-1}^n}^{t_i^n} \int_{t_{i-1}^n}^t  \frac{1}{(T-u)^{2\theta}} dudt.
\end{align*}
In the case $(L1)$ for our lower estimates \cite{Geiss15} Remark 6.6 implies that 
$$
\frac{1}{\sqrt{n}} \sum_{i\in I_n} \int_{t_{i-1}^n}^{t_i^n}\int_{t_{i-1}^n}^{t} \frac{1}{(T-u)^{2\theta}} du dt \leq \frac{1}{\sqrt{n}} \frac{C}{n^{1/2 + \eps}}
$$
for some $\eps > 0$. In the case $(L2)$, Lemma \ref{tthe:geiss} gives
$$
\frac{1}{\sqrt{n}} \sum_{i\in I_n} \int_{t_{i-1}^n}^{t_i^n}\int_{t_{i-1}^n}^{t} \frac{1}{(T-u)^{2\theta}} du dt \leq \frac{1}{\sqrt{n}} \frac{C}{n}. 
$$
The term containing $H^2$ can be bounded from below as
\begin{align*}
\liminf_{n\to \infty} n \sum_{i \in I_n} \int_{t_{i-1}^n}^{t_i^n} (t_i^n - t)\frac{H^2(t)}{4} dt 
&\geq \liminf_{n\to \infty} \frac{C_H}{4} n \sum_{i\in I_n} \frac{(t_i^n - t_{i-1}^n)^2}{2}\\
&\geq \liminf_{n\to \infty} \frac{C_H}{8} \vs \sum_{i\in I_n} (t_i^n - t_{i-1}^n)\os^2\\
&=\frac{C_H}{8}(B-A)^2.
\end{align*}
This proves the estimate.
\stopproof

\vspace*{2em}
\section{Examples}\label{cha:esim}
In this chapter we give two examples as an application of our results. For simplicity, we
consider a diffusion 
$$
X_t^i = x_0^i + \sum_{j=1}^d \int_0^t \sigma_{ij}(X_u) dW_u^j,\ i=1,\ldots,d
$$
in the case $(C2)$. Let $0=t_0 \leq t_1 \leq \cdots \leq t_n =T$ be any deterministic time-net on $[0,T]$. By $(t_i^\eta)_{i=0}^n$ we denote the time-net
$$
(t_i^\eta)_{i=0}^n = \vs T\vs 1 - \vs 1 - \frac{i}{n}\os^\frac{1}{1-\eta} \os\os_{i=0}^n \text{ and }
\vaalto \begin{array}{ll}
    \eta =0, & \theta \in [0,\frac{1}{2}) \\
    \eta\in (2\theta - 1, 1), & \theta \in [\frac{1}{2},1) \\
  \end{array}\right.
$$
where $\theta$ is from Theorem \ref{the:paatulos} equation (\ref{eq:tulos2oletus}).
\begin{exam}
For a European digital option with strike price $K > 0$,
$$
f(x):= \indi{\sum_{i=1}^d \lambda_i x_i\geq K}(x), \text{ where } \lambda_1,\ldots, \lambda_d > 0, 
$$
the approximation rate is $n^{-1/4}$ if equidistant time-nets are used \cite[Theorem 2.1]{Temam}.
Our main result Theorem \ref{the:paatulos} gives that this option can be approximated by the rate $n^{-1/2}$, more precisely 
$$
\vs \odotu \sup_{t\in[0,T]} \vits \sum_{i=1}^n \sum_{k=1}^d \int_{t_{i-1}^\eta \wedge t}^{t_i^\eta \wedge t} \vs\deri{x_k} F(u,X_u) - \deri{x_k} F(t_{i-1}^\eta,X_{t_{i-1}^\eta})\os dX_u^k \oits^2\os^\frac{1}{2} \leq \frac{D_1}{\sqrt{n}}
$$
and
$$
\frac{1}{D_2} \leq \liminf_{n\to\infty} \sqrt{n} \vs \odotu \vits \sum_{i=1}^n \sum_{k=1}^d \int_{t_{i-1}^\eta}^{t_i^\eta} \vs\deri{x_k} F(u,X_u) - \deri{x_k} F(t_{i-1}^\eta,X_{t_{i-1}^\eta})\os dX_u^k \oits^2\os^\frac{1}{2}
$$
for all $\eta \in (1/2, 1)$.
Assumption (\ref{eq:tulos2oletus}) follows from \cite[Proposition 3.3]{Temam} for $\theta = 3/4$ and (\ref{eq:oletus_H}) is due to \cite[equation (4.2)]{Temam}.
\end{exam}

\begin{exam}\label{exa:riippumattomuus}
In this example we will show how our main result Theorem \ref{the:paatulos} applies to the case of $\sigma$ being a $d\times d$ diagonal matrix with assumption $\sigma_{ii}(x) = \sigma_{ii}(x_i)$ and $d\geq 2$. As a special case of mixing different type of options we also apply Theorem \ref{the:paatulos} to the pay-off function
$$
f(x) := (x_1 - K_1)_+\ (x_2 - K_2)_+^\alpha\ \indi{[K_3,\infty)}(x_3)
$$
and show that for this pay-off one has the approximation rate $n^{-1/2}$ when approximation is optimized over deterministic time-nets of cardinality n+1.

Let us first turn to the case of $d\times d$ diagonal matrix.
If $d\geq 2$, $(\sigma_{ij})_{i,j=1}^d$ is a diagonal matrix, and if $\sigma_{ii}(x) = \sigma_{ii}(x_i)$, then the transition density of the process $Y$ can be written as the product of the transition densities of the process $Y^i$, i.e.
$$
\Gamma_Y(t,y,\xi) = \prod_{i=1}^d \Gamma_{Y^i}(t,y_i,\xi_i).
$$
Assume that $f(x) = \prod_{i=1}^d f_i(x_i)$, where the functions $f_i$ are of at most polynomial growth. The definition of $F$ implies that
$$
F(t,x) = \prod_{i=1}^d F_i(t,x_i)
$$
with
$$
F_i(t,x_i):= \int_\real \Gamma_{Y^i}(T-t, \log x_i, \xi_i) f_i(e^{\xi_i}) d\xi_i,
$$
and the second order derivatives of the function $F$ can be written as
$$
\deeri[2]{x_i x_j} F(t,x) = \vaalto
  \begin{array}{ll}
    \deri{x_i}F_i(t,x_i)\deri{x_j}F_j(t,x_j) \prod_{m=1,m\not=i, m\not=j}^d F_m(t,x_m) & i\not= j,  \\
     \deeri[2]{x_i x_i}F_i(t,x_i) \prod_{m=1,m\not=i}^d F_m(t,x_m) , & i=j. \\
  \end{array}\right.
$$
Assume that there exist $C>0$ and $\theta_i\in [0,1)$ for all $i=1,\ldots, d$ such that
\begin{equation}\label{eq:esim_tulo}
\odotu \vhaka Q_i^2(X_t) \deeri[2]{x_i x_i}F_i(t,X_t^i)\ohaka^2 \leq \frac{C}{(T-t)^{2\theta_i}}.
\end{equation}
Theorem \cite[Theorem 2.3]{Geissit} gives that 
\begin{equation}\label{eq:aste1_deriv_raj}
\sup_{t\in [0,T)} (T - t)^\delta \sqrt{t} \vs \odotu \vits \sigma_{ii}(X_t) \deri{x_i} F_i(t,X_t^i)\oits^2 \os^\frac{1}{2} <\infty
\end{equation}
if and only if
\begin{equation}\label{eq:aste2_deriv_raj}
\sup_{t\in [0,T)} (T - t)^\delta \vs \int_0^t \odotu \vits \sigma_{ii}^2(X_u) \deeri[2]{x_i x_i} F_i(u,X_u^i) \oits^2 du \os^\frac{1}{2} <\infty,
\end{equation}
where $\delta \in [0,1/2)$.
Now we assume, without loss of generality, that $\theta_i \in (1/2,1)$. For $\alpha = \beta$ we get that
\begin{align*}
&\odotu \vhaka A_{\alpha\alpha}(X_t) A_{\beta \beta}(X_t) \vs \deeri[2]{x_\alpha x_\beta} F(t,X_t)\os^2\ohaka\\
&=\odotu \vhaka \sigma_{\alpha\alpha}^2(X_t^\alpha)  \vits \deeri[2]{x_\alpha x_\alpha} F(t,X_t^\alpha)\oits^2\ohaka \prod_{m=1,m\not=\alpha}^d \odotu \vits F_m(t,X_t^m)\oits^2\\
&\leq \odotu \vhaka \sigma_{\alpha\alpha}^2(X_t^\alpha)  \vits \deeri[2]{x_\alpha x_\alpha} F(t,X_t^\alpha)\oits^2\ohaka \prod_{m=1,m\not=\alpha}^d \odotu \vits f_m(X_T^m)\oits^2
\end{align*}
which is at most of order $(T-t)^{2\theta_\alpha}$. For $\alpha \not= \beta$ we get that
\begin{align*}
&\odotu \vhaka A_{\alpha\alpha}(X_t) A_{\beta \beta}(X_t) \vs \deeri[2]{x_\alpha x_\beta} F(t,X_t)\os^2\ohaka\\
&\leq \odotu \vhaka \sigma_{\alpha\alpha}^2(X_t^\alpha)  \vits \deri{x_\alpha } F(t,X_t^\alpha)\oits^2\ohaka
 \odotu \vhaka \sigma_{\beta\beta}^2(X_t^\beta)  \vits \deri{x_\beta } F(t,X_t^\beta)\oits^2\ohaka \times \\
&\hspace*{2em}\times \prod_{m=1,m\not=\alpha,m\not=\beta}^d \odotu \vits f_m(X_T^m)\oits^2.
\end{align*}
The implication (\ref{eq:aste2_deriv_raj}) $\Rightarrow$ (\ref{eq:aste1_deriv_raj}) implies, for $\delta_i := (\theta_i - 1/2) \in (0,1/2)$, that
\begin{equation*}
\sup_{t\in [0,T)} (T - t)^{\delta_i} \sqrt{t} \vs \odotu \vits \sigma_{ii}(X_t) \deri{x_i} F_i(t,X_t)\oits^2 \os^\frac{1}{2} <\infty.
\end{equation*}
Using \cite[Lemma 5.2]{Geissit} one can remove the factor $\sqrt{t}$, so that 
\begin{equation*}
\sup_{t\in [0,T)} (T - t)^{2(\delta_\alpha + \delta_\beta)}  \odotu \vits A_{\alpha \alpha}(X_t) A_{\beta \beta}(X_t) \deeri[2]{x_\alpha x_\beta} F(t,X_t) \oits^2 <\infty.
\end{equation*}
Putting all estimates together, we find a $\theta \in [0,1)$ such that
\begin{equation*}
\odotu \vhaka A_{\alpha\alpha}(X_t) A_{\beta \beta}(X_t) \vs \deeri[2]{x_\alpha x_\beta} F(t,X_t)\os^2\ohaka \leq \frac{D}{(T-t)^{2\theta}}.
\end{equation*}
Looking at the above computations, one can take $\theta := \max\vaalto \theta_1,\ldots, \theta_d, 1/2\oaalto$ without the assumption $\theta_i \in (1/2, 1)$.

Let us go back to our $3$-dimensional example of mixing different type of options. Assume that $\sigma$ is a $3\times 3$ matrix defined by
$$
\sigma_{ij}(x) = \vaalto
  \begin{array}{ll}
    0, & i\not= j,  \\
    x_i, & i=j, 
  \end{array}\right.
$$
and that $x_0 = (1,1,1)$.
Define the pay-off function $f$ as above by
$$
f(x_1,x_2,x_3) = f_1(x_1)f_2(x_2)f_3(x_3):= (x_1 - K_1)_+\ (x_2 - K_2)_+^\alpha\ \indi{[K_3,\infty)}(x_3),
$$
where $K_i >0$, $i=1,2,3$ and $\alpha \in (0,\frac{1}{2})$. For $F_1$ one can compute
$$
\deeri[2]{x_1 x_1} F_1(t,x_1) = \frac{1}{x_1\sqrt{T-t}} \frac{1}{\sqrt{2\pi}} \exp\vhaka - \frac{\vs \frac{\log \vs\frac{x_1}{K_1}\os + \frac{T-t}{2}}{\sqrt{T-t}}\os^2}{2}\ohaka
$$
and
$$
\odotu\vhaka Q_1^2(X_t) \deeri[2]{x_1 x_1}F_1(t,X_t^1)\ohaka^2 = \frac{K_1}{2\pi \sqrt{T^2 - t^2}} \exp\vhaka -\frac{( T/2 + \log(K_1))^2}{T+t}\ohaka.
$$
This implies that one can choose $\theta_1=\frac{1}{4}$. For $F_2$, we can choose $\theta_2 = \frac{3-2\alpha}{4}$ and for $F_3$, $\theta_3 = \frac{3}{4}$ (cf. \cite[Lemma 1 and Lemma 2]{Gobet_Temam}). Now Theorem \ref{the:paatulos} gives that
$$
\vs \odotu \sup_{t\in[0,T]} \vits \sum_{i=1}^n \sum_{k=1}^d \int_{t_{i-1}^\eta \wedge t}^{t_i^\eta \wedge t} \vs\deri{x_k} F(u,X_u) - \deri{x_k} F(t_{i-1}^\eta,X_{t_{i-1}^\eta})\os dX_u^k \oits^2\os^\frac{1}{2} \leq \frac{D_1}{\sqrt{n}}
$$
for all $\eta \in (1/2, 1)$.
Under the assumptions of this example we have that
\begin{align*}
H^2(u) &= \odotu \sum_{\alpha, \beta = 1}^3 \vs \sigma_{\alpha \alpha}(X_u) \sigma_{\beta \beta} (X_u) \deeri[2]{x_\alpha x_\beta} F(u,X_u)\os^2\\
&\geq \odotu \sum_{\alpha = 1}^3 \vs \sigma_{\alpha \alpha}^2(X_u) \deeri[2]{x_\alpha x_\alpha} F(u,X_u)\os^2\\
&\geq \frac{1}{C_1} \sum_{\alpha = 1}^3 \odotu \vits X_u^\alpha\oits^4 \vits \deeri[2]{x_\alpha x_\alpha}F(u,X_u^\alpha)\oits^2 \prod_{m\not= \alpha}\odotu\vits F_m(u,X_u^m)\oits^2. 
\end{align*}
Since $f_1, f_2$ and $f_3$ are not almost surely linear and since
$$
u\mapsto \odotu \vits X_u^\alpha\oits^4 \vits \deeri[2]{x_\alpha x_\alpha}F_\alpha(u,X_u^\alpha)\oits^2
$$
is continuous and increasing, \cite[ proof of Proposition 2.1]{Geiss_BMO}, \cite[Theorem 4.6]{Geiss15} implies
$$
\sup_{u \in [0,T)} \odotu \vits X_u^\alpha\oits^4 \vits \deeri[2]{x_\alpha x_\alpha}F_\alpha(u,X_u^\alpha)\oits^2 >0. 
$$
Moreover,
$$
\lim_{u\nearrow T} \odotu \vs F_\alpha(u,X_u^\alpha)\os^2 = \odotu f_\alpha^2 (X_T^\alpha) > 0. 
$$
Hence Theorem \ref{the:paatulos} gives that
$$
\frac{1}{D_2} \leq \liminf_{n\to\infty} \sqrt{n} \vs \odotu \sup_{t \in [0,T]} \vits \sum_{i=1}^n \sum_{k=1}^d \int_{t_{i-1}^{n,\eta} \wedge t}^{t_i^{n,\eta} \wedge t} \vs\deri{x_k} F(u,X_u) - \deri{x_k} F(t_{i-1}^{n,\eta},X_{t_{i-1}^{n,\eta}})\os dX_u^k \oits^2\os^\frac{1}{2}.
$$
If we take $f(x_1,x_2) = f_1(x_1)f_2(x_2)$, then we can choose $\theta = \theta_2 < 3/4$ and get
$$
\frac{1}{D_2} \leq \liminf_{n\to\infty} \sqrt{n} \vs \odotu \sup_{t \in [0,T]} \vits \sum_{i=1}^n \sum_{k=1}^d \int_{t_{i-1}^n\wedge t}^{t_i^n\wedge t} \vs\deri{x_k} F(u,X_u) - \deri{x_k} F(t_{i-1},X_{t_{i-1}})\os dX_u^k \oits^2\os^\frac{1}{2},
$$
for any sequence of time-nets with $\sup_{i=1,\ldots,n}(t_i^n - t_{i-1}^n) \leq C/n$.
\end{exam}

\vspace*{2em}
\section*{Appendix}
\setcounter{section}{5}
\setcounter{lemma}{0}
\setcounter{equation}{0}
\begin{theorem}[Theorem 8. p. 263 in \cite{Friedman64}, Theorem 5.4. p. 149 in \cite{Friedman75}]\label{the:siirtymatod}
For $\hat{b},\ \hat{\sigma}$ with $\hat{\sigma}\hat{\sigma}^T$ uniformly elliptic, there exists a transition density $\Gamma:(0,T]\times \real^d \times \real^d \to [0, \infty) \in C^\infty$ such that $\tn(Y_t \in B) = \int_B \Gamma(t,y,\xi)d\xi$, for $t \in (0,T]$ and $B\in \borel{\real^d}$, where $Y = (Y_t)_{t\in [0,T]}$ is the strong solution of the SDE (\ref{eq:SDE_Y}) starting from $y$: Moreover, the following are satisfied:
\begin{itemize}
\item[(i)] For $(s,y,\xi) \in (0,T]\times \real^d \times \real^d$ one has
 $$
 \deri{s} \Gamma(s,y,\xi) = \frac{1}{2} \sum_{k,l=1}^d \sum_{j=1}^d \hat{\sigma}_{kj}(y) \hat{\sigma}_{lj}(y) \deeri[2]{y_k y_l} \Gamma(s,y,\xi) + \sum_{i=1}^d \hat{b}_i(y) \deri{y_i}\Gamma(s,y,\xi).
 $$
 \item[(ii)] For $a \in \vaalto 0,1,2,\ldots\oaalto$ and multi-indices $b$ and $c$ there exist positive constants $C$ and $D$, depending only on $a,b,c$ and $d$, such that
 $$
 \vits \frac{\partial^{a+|b|+|c|}}{\partial^a t\ \partial^b y\ \partial^c \xi} \Gamma(t,y,\xi) \oits \leq \frac{C}{t^{(d+2a + |b| + |c|)/2}} e^{-D \frac{\vnorm y-\xi \onorm^2}{t}},
 $$
\end{itemize}
where $\vnorm \cdot \onorm$ is the Euclidean norm. 
\end{theorem}

\begin{remark}
In references \cite{Friedman64} and \cite{Friedman75} it is assumed that functions $\hat{b}$ and $\hat{\sigma}$ are Lipschitz continuous. Here it follows from the assumption (\ref{eq:kerroinfunktio_oletus}).
\end{remark}

\begin{theorem}[Gronwall's Lemma, \cite{yor}]\label{the:gronwall}
If, for $ t_0\leq t\leq t_1$, $ \phi(t)\geq 0$ is a continuous function such that
$$
 \phi(t)\leq K+L\int_{t_0}^t \phi(s)ds 
$$
for $t_0\leq t\leq t_1$ where $K,L\geq 0$, then
$$
 \phi(t)\leq K e^{L (t-t_0)} 
$$
on $ t_0\leq t\leq t_1$.
\end{theorem}

\begin{lemma}\label{lemma:pysaytys}
Let $0\leq a < b < T$ and define
$$
\phi_{kl}(u,x) := \vs \deri{x_k} F(u,x) - v_a^k\os \sigma_{kl}(x),\ u \in [0,T),\ x \in E,
$$
where $v_a^k$ is an $\saF_a$-measurable random variable and assume that
$$
\odotu \sup_{u\in [a,b]}\vhaka \phi_{kl}^2(u,X_u) + \vits \vs \oA
\phi_{kl}^2\os (u,X_u)\oits +  \sum_{m=1}^d \vits \deri{x_m}\phi_{kl}^2(u,X_u) b_m(X_u) \oits \ohaka < \infty.
$$
Then for $s\in [a,b]$ one has
\begin{align*}
\odotu \phi_{kl}^2 (s,X_s) = &\odotu \phi_{kl}^2(a,X_a) + \int_a^s \odotu \vs\oA  \phi_{kl}^2\os (u,X_u) du + \\
\hspace*{2em}& + \int_a^s \odotu \sum_{m=1}^d \vs \deri{x_m}\phi_{kl}^2(u,X_u)\os b_m(X_u) du.
\end{align*}
\end{lemma}

%\proof
\myproof
By It\^o's formula we obtain
\begin{align*}
 \phi_{kl}^2 (s,X_s) &= \phi_{kl}^2(a,X_a) + \int_a^s  \vs\oA  \phi_{kl}^2\os (u,X_u) du\\
&+ \sum_{m=1}^d \int_a^s \deri{x_m} \phi_{kl}^2 (u,X_u) dX_u^m.
\end{align*}
Define
$$
S_n^m := \inf \vaalto r\in [a,s]\vits \int_a^r\sum_{j=1}^d\vhaka
\deri{x_m}\phi_{kl}^2(u,X_u)\ohaka^2 \sigma_{mj}^2(X_u)du > n\right. \oaalto \wedge
s
$$
and
$$
S_n := \min\vaalto S_n^m, m\in \vaalto 1,\ldots, d\oaalto \oaalto.
$$
This implies that
$$
 \int_a^{S_n}\sum_{j=1}^d\vhaka
\deri{x_m}\phi_{kl}^2(u,X_u)\ohaka^2 \sigma_{mj}^2(X_u)du \leq n
$$
and
$$
\odotu \int_a^{S_n} \deri{x_m} \phi_{kl}^2(u,X_u) \sigma_{mj}(X_u) dW_u^j = 0,
$$
for $n \in \natu$, $m\in \vaalto 1, \ldots, d\oaalto$ and $j\in \vaalto 1, \ldots, d\oaalto$.
Dominated convergence gives
\begin{align*}
\odotu \phi_{kl}^2(s,X_s) &= \lim_{n\to \infty} \odotu \phi_{kl}^2(S_n,X_{S_n})\\
&= \lim_{n\to \infty} \vhaka \odotu \phi_{kl}^2(a,X_a) + \odotu \int_a^{S_n}  \vs\oA  \phi_{kl}^2\os (u,X_u)
du +\right.\\
&\hspace*{4em} +\left. \odotu \sum_{m=1}^d \int_a^{S_n} \deri{x_m} \phi_{kl}^2 (u,X_u)
dX_u^m\ohaka\\
&= \odotu \phi_{kl}^2(a,X_a) + \odotu \int_a^s \vs\oA  \phi_{kl}^2\os
(u,X_u)du+\\
&\hspace*{2em} + \odotu \sum_{m=1}^d\int_a^s \deri{x_m} \phi_{kl}^2 (u,X_u) b_m(X_u)du.
\end{align*}
\stopproof

\begin{lemma}\label{lemma:rajoittuneisuus}
If for all $x \in E$
$$
\vits \frac{\partial^s}{\partial_{x_\beta}^q \partial_{x_\alpha}^r} \sigma_{ij}(x)\oits \leq C \frac{Q_i(x)}{Q_\beta^q(x) Q_\alpha^r(x)}, \ q+r = s,\ q,r,s\in \vaalto 0,1,2\oaalto,
$$
for some $C > 0$,
then for all $0\leq a \leq b < T$ and $k,l \in \vaalto 1,\ldots, d\oaalto$ we have that
\begin{equation*}
\odotu \sup_{u\in [a,b]} \phi_{kl}^2(u,X_u) < \infty,
\end{equation*}
\begin{equation*}
\odotu \sup_{u\in [a,b]} \vits Q_m(X_u) \deri{x_m} \phi_{kl}(u,X_u)\oits^2  < \infty, \ m = 1,\ldots, d,
\end{equation*}
\begin{equation*}
\odotu \sup_{u\in [a,b]} \vits Q_m(X_u)  \deri{x_m} \phi_{kl}^2(u,X_u) \oits < \infty, \ m = 1,\ldots, d,
\end{equation*}
and 
\begin{equation*}
\odotu \sup_{u\in [a,b]} \vits \vs \oA \phi_{kl}^2\os (u,X_u)\oits < \infty,
\end{equation*}
where
$$
\phi_{kl}(u,x) = \vs \deri{x_k}F(u,x) - \deri{x_k}F(a,X_a)\os \sigma_{kl}(x),\ u\in [a,b]. 
$$
\end{lemma}

\myproof
This proof uses the same notation for different constants.
Equation (\ref{eq:ey_F_1_deri}) implies that the random variable $\phi_{kl}^2(u,X_u)$ can be bounded by
\begin{align*}
\phi_{kl}^2(u,X_u) &= \vs \deri{x_k}F(u,X_u) - \deri{x_k} F(a,X_a)\os^2 \sigma_{kl}^2(X_u)\\
&\leq C \vs \deri{x_k}F(u,X_u) - \deri{x_k} F(a,X_a)\os^2 Q_k^2(X_u)\\
&\leq C \vs \vits Q_k(X_u) \deri{x_k}F(u,X_u)\oits^2 + \vits Q_k(X_a) \deri{x_k}F(a,X_a)\oits^2 \frac{Q_k^2(X_u)}{Q_k^2(X_a)}\os\\
&\leq C \vhaka\sup_{u' \in [a,b]} \vs 1 + \vnorm X_{u'}\onorm^q\os^2\ohaka \vs 1 + \frac{Q_k^2(X_u)}{Q_k^2(X_a)}\os.
\end{align*}
Applying H\"older's inequality we get that
\begin{align*}
\odotu \sup_{u \in [a,b]}\phi_{kl}^2(u,X_u) 
&\leq C \odotu \vhaka \sup_{u \in [a,b]} \vs 1 + \vnorm X_u\onorm^q\os^2  \sup_{u \in [a,b]}\vs 1 + \frac{Q_k^2(X_u)}{Q_k^2(X_a)}\os\ohaka\\
&\leq C  \vs \odotu  \sup_{u \in [a,b]} \vs 1 + \vnorm X_u\onorm^q\os^4\os^\frac{1}{2} \vs \odotu  \sup_{u \in [a,b]}\vs 1 + \frac{Q_k^2(X_u)}{Q_k^2(X_a)}\os^2 \os^\frac{1}{2}.
\end{align*}
Equation (\ref{eq:X_normi_rajoitettu}) gives that $\odotu  \sup_{u \in [a,b]} \vs 1 + \vnorm X_u\onorm^q\os^4$ is finite. 
In the case $(C1)$ it is trivial that the latter term is finite. Let us now turn to the case $(C2)$. Theorem \ref{the:siirtymatod} implies that 
\begin{align}\label{eq:odotu_X_-p_rajoitettu}
\odotu (X_u^k)^{-p} &= \odotu e^{-p Y_u^k}\\\nonumber
&= \int_{\real^d} e^{-p y_k} \Gamma_{Y}(u,y_0, y)dy\\\nonumber
&\leq C\  \exp\vhaka -p y_0^k + \frac{1}{2} \frac{u}{2D} p^2\ohaka < \infty
\end{align}
for all $p \in [0,\infty)$ and some $C>0$ and $D>0$. H\"older's inequality now implies that
\begin{align*}
\odotu  \sup_{u \in [a,b]}\vs 1 + \frac{Q_k^2(X_u)}{Q_k^2(X_a)}\os^2
&\leq 2  \vs 1 + \odotu \sup_{u \in [a,b]}\frac{Q_k^4(X_u)}{Q_k^4(X_a)}\os\\
&\leq 2  \vs 1 + \vs \odotu \sup_{u \in [a,b]}(X_u^k)^8\os^\frac{1}{2} \vs \odotu \frac{1}{(X_a^k)^8}\os^\frac{1}{2}\os,
\end{align*}
which is finite by (\ref{eq:X_normi_rajoitettu}) and (\ref{eq:odotu_X_-p_rajoitettu}).

Straightforward calculation gives that
\begin{align*}
&\odotu \sup_{u\in [a,b]} \vits Q_m(X_u) \vs \deri{x_m}\phi_{kl}\os(u,X_u) \oits^2\\
&= \odotu \sup_{u\in [a,b]} \vhaka Q_m^2(X_u) \vits\deeri[2]{x_m x_k}F(u,X_u) \sigma_{kl}(X_u) +\right.\right.\\
&\hspace*{8em} \left. \left.+ \vs \deri{x_k}F(u,X_u) - \deri{x_k}F(a,X_a)\os \vs\deri{x_m}\sigma_{kl}\os(X_u) \oits^2\ohaka\\
&\leq C \vs \odotu \sup_{u\in [a,b]} Q_m^2(X_u) Q_k^2(X_u) \vits\deeri[2]{x_m x_k}F(u,X_u)\oits^2+\right.\\
&\hspace*{4em} \left. + \odotu \sup_{u\in [a,b]} Q_k^2(X_u) \vits \deri{x_k}F(u,X_u) - \deri{x_k}F(a,X_a)\oits^2\os,
\end{align*}
and this is finite by equations (\ref{eq:X_normi_rajoitettu}) and (\ref{eq:ey_F_2_deri}) and the above argument.

H\"older's inequality together with the above gives that
\begin{align*}
&\odotu \sup_{u\in [a,b]} \vits Q_m(X_u) \vs\deri{x_m}\phi_{kl}^2\os(u,X_u) \oits\\
&= 2 \odotu \sup_{u\in [a,b]} \vits Q_m(X_u)  \phi_{kl} (u,X_u) \vs \deri{x_m}\phi_{kl} \os (u,X_u)\oits\\
&\leq 2\vs\odotu \sup_{u\in [a,b]} \vits \phi_{kl} (u,X_u)\oits^2\os^\frac{1}{2} \vs\odotu \sup_{u\in [a,b]}  Q_m^2(X_u)\vits \vs \deri{x_m}\phi_{kl}\os (u,X_u) \oits^2\os^\frac{1}{2} < \infty.
\end{align*}

For the last part of the proof, equation (\ref{eq:A_phi_toiseen}) gives that
$$
\odotu \sup_{u\in [a,b]} \vits \vs \oA \phi_{kl}^2\os(u,X_u)\oits < \infty
$$
if
$$
\odotu \sup_{u\in [a,b]}  \vs \oA \phi_{kl}\os^2(u,X_u) < \infty.
$$
This follows from equation (\ref{eq:equ2}) and the above arguments.
\stopproof

\begin{lemma}\label{tthe:geiss}
Assume that a Borel-measurable function $\vph:[0,T) \to [0,\infty)$ satisfies
$$
\vph(u) \leq \frac{C}{(T-u)^\theta}, \ u\in [0,T),
$$
for some $C >0$ and some $\theta \in [0,1)$.
Then there exists a constant $C'>0$ such that 
$$
\sum_{t_i \in \Tau_n^\eta} \int_{t_{i-1}^\eta}^{t_i^\eta} \int_{t_{i-1}^\eta}^u \vph^2(s) ds du \leq \frac{C'}{n},
$$
where 
$$
\tau_n^\eta  = (t_i^\eta)_{i=0}^n := \vs T\vs 1 - \vs 1 - \frac{i}{n}\os^\frac{1}{1-\eta} \os\os_{i=0}^n \text{ and }
\vaalto \begin{array}{ll}
    \eta =0, & \theta \in [0,\frac{1}{2}) \\
    \eta\in (2\theta - 1, 1), & \theta \in [\frac{1}{2},1). \\
  \end{array}\right.
$$
\end{lemma}

\myproof
Lemma follows from \cite[Lemma 4.14 and Proposition 4.16]{Gei15pre}. 
\stopproof


\vspace*{2em}






\begin{thebibliography}{1}

\bibitem{Friedman64}
Friedman, A. (1964):
\newblock \emph{Partial Differential Equations of Parabolic Type}.
\newblock Prentice-Hall.

\bibitem{Friedman75}
Friedman, A. (1975):
\newblock \emph{Stochastic Differential Equations and Applications}.
\newblock Vol. 1, New York: Academic Press


\bibitem{Geissit}
Geiss, C., Geiss, S. (2004):
\newblock On approximation of a class stochastic integrals  and interpolation,
\newblock {\em Stochastics and Stochastic Reports}, Vol. 76, No. 4, 339-362.


\bibitem{Gei15pre}
Geiss, S. (1999):
\newblock On quantitative approximation of stochastic integrals with respect to
  the geometric brownian motion,
\newblock {\em Report Series: SFB Adaptive Information Systems and Modelling in
  Econo\-mics and Management Science, Vienna University}, 43.

\bibitem{Geiss15}
Geiss, S. (2002):
\newblock Quantitative approximation of certain stochastic integrals,
\newblock {\em Stochastics and Stochastic Reports}, 73, 241-270.

\bibitem{Geiss_BMO}
Geiss, S. (2002):
\newblock On the approximation of stochastic integrals and weighted BMO,
\newblock {\em Stochastic Processes and Related Topics}.

\bibitem{Geiss_Hujo}
Geiss, S., Hujo, M. (2007):
\newblock Interpolation and approximation in $L_2(\gamma)$
\newblock {\em Journal of Approximation Theory} 144, 213 - 232.



\bibitem{Gobet_Temam}
Gobet, E., Temam, E. (2001):
\newblock Discrete time hedging errors for options with irregular payoffs,
\newblock {\em Finance and Stochastics}, 5, 357-367.

\bibitem{Kar-Shreve1}
Karatzas, I., Shreve, S. (1988):
\newblock \emph{Brownian Motion and Stochastic Calculus}.
\newblock Springer.

\bibitem{Nualart}
Nualart, D. (1995):
\newblock  \emph{The Malliavin calculus and related topics}.
\newblock Springer-Verlag.


\bibitem{yor}
Revuz, D., Yor, M. (1991): \newblock \emph{Continuous martingales and Brownian
  motion}. Springer-Verlag.


\bibitem{Temam}
Temam, E. (2003):
\newblock Analysis of error with malliavin calculus: Application to hedging,
\newblock {\em Mathematical Finance}, 13 , 1, 201-214.


\bibitem{ZhangR}
Zhang, R. (1999):
\newblock \emph{Couverture approch\'ee des options {E}urop\'eennes}.
\newblock PhD thesis, Ecole Nationale des Ponts et Chauss\'ees.



\end{thebibliography}
\end{document}